\documentclass[12pt]{amsart}

\usepackage[arrow,dvips,matrix,arrow,ps,color,line,curve,frame]{xy}

\usepackage[arrow,matrix]{xy}
\usepackage{amssymb}

\advance\textheight1cm

\let\oldlabel=\label
\def\prellabel{\marginparsep=1em
    \def\label##1{\oldlabel{##1}\ifmmode\else\ifinner\else
         \marginpar{{\footnotesize\ \\ \tt
                    ##1}}\fi\fi}}

\advance\headheight1.15pt

\def\supp{{\operatorname{supp}}}
\def\pyr{{\operatorname{pyr}}}
\def\cc{{\operatorname{\mathfrak c}}}
\def\inte{{\operatorname{int}}}
\def\gp{{\operatorname{gp}}}
\def\Ker{{\operatorname{Ker}}}

\def\SL{\operatorname{SL}}
\def\conv{\operatorname{conv}}
\def\E{\operatorname{E}}
\def\st{\operatorname{St}}
\def\rank{\operatorname{rank}}
\def\Ker{{\operatorname{Ker}}}

\def\vert{\operatorname{vert}}
\def\Im{\operatorname{Im}}
\def\GL{\operatorname{GL}}

\def\conv{\operatorname{conv}}
\def\Q{\qedsymbol\kern1pt}
\def\sqq#1#2{{\hbox{\kern 0.5pt\vbox{\hbox{\kern#1ex\vrule width #2pt height #1ex}%
     \hrule height #2pt}\kern 1.5pt}}}

\def\vt{{\kern1pt|}}
\def\RR{{\mathbb R}}

\def\QQ{{\mathbb Q}}
\def\ZZ{{\mathbb Z}}
\def\NN{{\mathbb N}}

\def\1{\text{\bf 1}}
\def\A{{\mathcal A}}
\def\B{{\mathcal B}}
\def\D{{\mathcal D}}
\def\M{{\mathcal M}}
\def\H{{\mathcal H}}
\def\G{{\mathcal G}}
\def\k{{\mathbf k}}

\let\epsilon=\varepsilon
\let\phi=\varphi
\let\theta=\vartheta

\newtheorem{lemma}{Lemma}[section]
\newtheorem{corollary}[lemma]{Corollary}
\newtheorem{theorem}[lemma]{Theorem}
\newtheorem{proposition}[lemma]{Proposition}

\theoremstyle{definition}

\newtheorem{remark}[lemma]{Remark}

\begin{document}

\title[Steinberg group of a monoid ring]
{The Steinberg group of a monoid ring,\\ nilpotence, and
algorithms}

\author{Joseph Gubeladze}

\address{Department of Mathematics,
San Francisco State university, San Francisco, CA 94132, USA}
\email{soso@math.sfsu.edu}

\subjclass[2000]{14M25, 19B14, 19C09, 20G35, 52B20}

\begin{abstract}
For a regular ring $R$ and an affine monoid $M$ the homotheties of
$M$ act nilpotently on the Milnor unstable groups of $R[M]$. This
strengthens the $K_2$ part of the main result of \cite{G5} in two
ways: the coefficient field of characteristic 0 is extended to any
regular ring and the stable $K_2$-group is substituted by the
unstable ones. The proof is based on a polyhedral/combinatorial
technique, computations in Steinberg groups, and a substantially
corrected version of an old result on elementary matrices by
Mushkudiani \cite{Mu}. A similar stronger nilpotence result for
$K_1$ and algorithmic consequences for factorization of high
Frobenius powers of invertible matrices are also derived.
\end{abstract}

\maketitle

\section{Introduction}\label{INTR}

\subsection{Main result} In the recent work \cite{G5} we proved
the following result. Let $\k$ be a field of characteristic 0, $M$
be an additive submonoid of $\ZZ^n$ without nontrivial units, and
$i$ be a nonnegative integer. Then for any element $x\in
K_i(\k[M])$ and any natural number $c\ge2$ there exists an integer
$j_x\ge0$ such that $(c^j)_*(x)\in K_i(\k)$ for all $j\ge j_x$.

\smallskip Here for a natural number $c$ the group endomorphism
 of $K_i(\k[M])$, induced by the monoid
endomorphism $M\to M$, $m\mapsto m^c$, is denoted by $c_*$.

\smallskip The motivation for this result is that it is a
natural higher version of the triviality of algebraic vector
bundles on affine toric varieties \cite{G1}, contains Quillen's
fundamental result on homotopy invariance, and easily extends to
global toric varieties. See the introduction of \cite{G5} for the
details.

\smallskip This result confirms the
\emph{nilpotence conjecture} for a special class of coefficients
rings. The conjecture asserts the similar nilpotence property of
higher $K$-groups of monoid algebras over \emph{any (commutative)
regular} coefficient ring.

\smallskip The main result in this paper is a stronger
unstable version of the nilpotence property for the functors
$K_{1,r}$ and $K_{2,r}$ for any regular coefficient ring.
Moreover, when the coefficient ring is a field the argument leads
to an algorithm for factorization of high `Frobenius powers' of
invertible matrices into elementary ones.

\smallskip In the special case of the polynomial rings
$\k[\ZZ_+^n]=\k[t_1,\ldots,t_n]$ the algorithmic study of
factorizations of invertible matrices has applications in signal
processing \cite{LiXW,PW}. The starting point here is Suslin's
well known paper \cite{Su}. In this special case there is no need
to take Frobenius powers of invertible matrices. However, a
$K$-theoretical obstruction shows that this is no longer possible
once we leave the class of free monoids, see Remark \ref{Cohn}.
Therefore, our algorithmic factorization is an optimal `sparse
version' of the existing algorithm for polynomial rings.

\smallskip Here is the main result:

\begin{theorem}\label{MAIN}
Let $M$ be a commutative cancellative torsion free monoid without
nontrivial units, $c\ge2$ a natural number, $R$ a commutative
regular ring and $\k$ a field. Then:
\begin{itemize} \item[(a)] For any element $z\in K_{2,r}(R[M])$,
$r\ge\max(5,\dim R+3)$, there exists an integer $j_z\ge0$ such
that
$$
(c^j)_*(z)\in K_{2,r}(R)=K_2(R),\qquad j\ge j_z.
$$
\item[(b)] For any matrix $A\in\GL_r(R[M])$, $r\ge\max(3,\dim
R+2)$, there exists an integer number $j_A\ge0$ such that
$$
(c^j)_*(A)\in\E_r(R[M])\GL_r(R),\qquad j\ge j_A.
$$
\item[(c)] There is an algorithm which for any matrix
$A=\SL_r(\k[M])$, $r\ge3$, finds an integer number $j_A\ge0$ and a
factorization of the form:
\begin{align*}
(c^{j_A})_*(A)=\prod_ke_{p_kq_k}(\lambda_k),\qquad
e_{p_kq_k}(\lambda_k)\in\E_r(\k[M]).
\end{align*}
\end{itemize}
\end{theorem}
Here:

\smallskip$\centerdot$ for a commutative ring $\Lambda$ its Krull dimension
is denoted by $\dim\Lambda$,

\smallskip$\centerdot$ $K_{2,r}(-)$ refers to the Milnor's $r$th unstable $K_2$,

\smallskip$\centerdot$ for a natural number $c$ the group endomorphisms
$\GL_r(R[M])\to\GL_r(R[M])$ and $K_{2,r}(R[M])\to K_{2,r}(R[M])$,
induced by the monoid endomorphism $M\to M$, $m\mapsto m^c$, are
both denoted by $c_*$.

\smallskip$\centerdot$ for two subgroup $H_1$ and $H_2$ of a group
$G$ we use the notation $H_1H_2=\{h_1h_2\ |\ h_1\in H_1,\ h_2\in
H_2\}$.

\begin{remark}\label{norealgorithm}
We do not give a detailed description of the actual algorithm,
mentioned in Theorem \ref{MAIN}(c). Rather, throughout the text,
we highlight the explicit nature of the proof of Theorem
\ref{MAIN}(b) which implies the possibility of converting the
argument into an implemented algorithm when the coefficients are
in a field.
\end{remark}

\begin{remark}\label{Kregular}
It is not difficult to show that the proof of Theorem \ref{MAIN},
given below, works for a more general class of rings of
coefficients. In fact, all one needs from the ring $R$ is the
validity of the claims (a), (b) and (c) for the polynomial
extension $R[t_1,\ldots,t_n]$, $n=\rank M$ -- a classically known
fact when $R$ is a regular ring, see Section \ref{Kbaground}.
\end{remark}

\begin{remark}\label{uniform}
We do not know whether there is a uniform bound $j_i$, depending
on $M$, $R$ and $i$, but not on $x\in K_i(R[M])$, such that
$(c^{j_i})_*(x)\in K_i(R)$. Nontrivial examples in \cite{G3}
indicate that such bounds may in fact exist, at least for $K_1$.
\end{remark}

\smallskip A word is in order on the previous results and the
proof of Theorem \ref{MAIN}.

\smallskip The proof of the nilpotence of
$K_i(\k[M])$ as given in \cite{G5} -- even in the case of Milnor's
$K_2$ -- uses a series of deep facts in higher $K$-theory of
rings, obtained from the early 1990s on (the most recent of which
is \cite{Cor}). The proof of Theorem \ref{MAIN}, given below,
makes no use of any of these results. It is based on computations
in $\E_r(R[M])$, essentially  due to Mushkudiani \cite{Mu}, and
similar computations in $\st_r(R[M])$. The explicit nature of
these computations is also the source of the algorithmic
consequences for $\SL_r(\k[M])$. Obviously, no such a pure
algebraic approach is possible for higher $K$-groups.

\smallskip Actually, the weaker stable version of Theorem
\ref{MAIN}(a) for $K_2$ is claimed in \cite{Mu} and the present
work grew up from our attempts to understand Mushkudiani's
argument. Eventually, what survived from \cite{Mu} is his
preliminary computations in the group of elementary matrices -- an
important technical fact whose corrected and stronger unstable
version is given in the last Section \ref{Mushk}; see Remarks
\ref{whycompl}, \ref{mushkwrong} and
\ref{mushksteparation}.\footnote{We also greatly simplify the
notation in \cite{Mu} -- already a challenge on its own right.}
The rest of the paper is devoted to the reduction of Theorem
\ref{MAIN} to this technical fact.

\smallskip In the course of the proof we also develop an
effective/algorithmic excision technique for the unstable $K_1$
and $K_2$-groups of monoid rings (Section \ref{redint}). It allows
us to circumvent Suslin-Wodzicki's excision theorem \cite{SuW} --
a result which is applicable only to stable groups and which was
essential in \cite{G5}.

\smallskip Finally, a comment on the result on $K_1$: the weaker
stable analog of Theorem \ref{MAIN}(b) is obtained in \cite{G2},
where we originally conjectured the nilpotence of the higher
$K$-theory of $R[M]$. But the essential difference between the two
approaches is that in the present paper we never invoke
\emph{Quillen's local-global patching}, \emph{Karoubi squares} and
\emph{Horrock's localizations at monic polynomials}, heavily used
in \cite{G1,G2,G5}. On the other hand, it should be mentioned that
the technique developed in \cite{G2} is crucial in the proof of
the nilpotence result for higher $K$-groups, see \cite[\S9]{G4}.

\subsection{Organization of the paper} To make the exposition as
self-contained as possible, the necessary $K$-theoretical
background, together with a further motivation for the main
result, is provided in Section \ref{Kbaground}. In Section
\ref{MonCon} we give a quick summary of the polyhedral approach to
commutative, cancellative, torsion free monoids, developed in our
study of $K$-theory of monoid rings.  An effective excision
technique for unstable $K_1$- and $K_2$-groups of monoid rings is
developed in Section \ref{redint}. In Section \ref{proofmain} we
introduce an inductive process, \emph{pyramidal descent}, on which
the proof of Theorem \ref{MAIN} is based. The main technical facts
that make this inductive process work, Theorems \ref{elparation}
and \ref{steparation}, are stated in Section \ref{ALMOST}. There
we also explain how \ref{steparation} follows from
\ref{elparation}. In Section \ref{PYRAMID} we show the validity of
pyramidal descent in the situation of Theorem \ref{MAIN}. Section
\ref{Mushk} presents a corrected version of Mushkudiani's proof of
Theorem \ref{elparation}.

\

\noindent{\emph{Acknowledgment.} I am grateful to the referee for
the thorough study of the paper and spotting various inaccuracies
in the original version, especially in Section \ref{Mushk}.

\section{$K$-theoretical background}\label{Kbaground}
Let $\Lambda$ be a ring and $r\ge2$ a natural number. For a pair
of natural numbers $1\le p,q\le r$ and an element
$\lambda\in\Lambda$ the  matrix with $\lambda$ on the
$pq$-position and 0s elsewhere will be denoted $a_{pq}(\lambda)$.

\smallskip \emph{The standard elementary matrices} over $\Lambda$
of order $r$ are defined as follows
\begin{align*}
e_{pq}(\lambda)=\1+a_{pq}(\lambda),\quad1\le p,q\le r,\quad
p\not=q,\quad\lambda\in \Lambda,
\end{align*}
where $\1$ is the unit matrix.

\smallskip The standard elementary matrices generate \emph{the
subgroup of elementary matrices} $\E_r(\Lambda)$ inside the
general linear group $\GL_r(\Lambda)$ of order $r$.

\smallskip Starting from now on all our rings are assumed to be
\emph{commutative}.

\smallskip It is known that $\E_r(\Lambda)\subset\GL_r(\Lambda)$ is
a normal subgroup as soon as $r\ge3$ \cite{Su}.

\smallskip The \emph{special linear group $\SL_r(\Lambda)$ of order $r$} is
defined to be the subgroup of $\GL_r(\Lambda)$ of the matrices
with determinant 1. Thus
$\E_r(\Lambda)\subset\SL_r(\Lambda)\subset\GL_r(\Lambda)$.

\smallskip Let $G_r$ denote any of the groups $\E_r(\Lambda)$,
$\SL_r(\Lambda)$, $\GL_r(\Lambda)$. \emph{The stable group $G$} is
defined to be the inductive limit of the diagram of groups
\begin{align*}
&G_2(\Lambda)\to\cdots\to G_r(\Lambda)\to G_{r+1}(\Lambda)\to\cdots,\\
&A\mapsto
\begin{pmatrix}
A&0\\
0&1
\end{pmatrix},\qquad A\in G_r(\Lambda).
\end{align*}
The \emph{Whitehead Lemma} says that
$\E(\Lambda)=\left[\GL(\Lambda),\GL(\Lambda)\right]$ \cite[Lemma
3.1]{Mi}. \emph{The Bass-Whitehead group} $K_1(\Lambda)$ is
defined by
$$
K_1(\Lambda)=\GL(\Lambda)/\E(\Lambda)=\GL(\Lambda)_{\text{ab}}=H_1(\GL(\Lambda),\ZZ).
$$
Its unstable versions are given by
$K_{1,r}(\Lambda)=\GL_r(\Lambda)/\E_r(\Lambda)$, $r\ge3$.

\smallskip The standard elementary matrices satisfy \emph{the
Steinberg relations}:
\begin{align*}
&\ e_{pq}(\lambda)\cdot e_{pq}(\mu)=e_{pq}(\lambda+\mu),\\
&[e_{pq}(\lambda),e_{qu}(\mu)]=e_{pu}(\lambda\mu),\quad p\not=u\\
&[e_{pq}(\lambda),e_{uv}(\mu)]=1,\quad p\not=v,\ q\not=u.
\end{align*}

\smallskip \emph{The unstable Steinberg group $\st_r(\Lambda)$} (over
$\Lambda$) is defined by the generators $x_{pq}(\lambda)$,  $1\le
p,q\le r$, $p\not=q$ and $\lambda\in \Lambda$, subject to the
corresponding Steinberg relations. The stable group $\st(\Lambda)$
is the inductive limit of the diagram
$\st_2(\Lambda)\to\st_3(\Lambda)\to\cdots$.

\smallskip \emph{The  Milnor $r$-th unstable group}
$K_{2,r}(\Lambda)$ is defined as the kernel of the canonical
surjective group homomorphism $\st_r(\Lambda)\to\E_r(\Lambda)$.
Passing to the inductive limits we get the short exact sequence of
the corresponding stable groups:
$$
1\to K_2(\Lambda)\to\st(\Lambda)\to\E(\Lambda)\to1.
$$
This is the sequence of a universal central extension of the
perfect group $\E(\Lambda)$ \cite[Theorem 5.10]{Mi}. Consequently,
$K_2(\Lambda)=H_2(\E(R),\ZZ)$.

Van der Kallen has shown \cite{K2} that the extension
$$
1\to K_{2,r}(\Lambda)\to\st_r(\Lambda)\to\E_r(\Lambda)\to1.
$$
is also universal central if $r\ge5$.

\smallskip All groups mentioned above, stable or unstable,
depend functorially on the underlying ring $\Lambda$.

\begin{theorem}\label{BHS}
Let $R$ be a regular ring. Then $K_i(R)=K_i(R[t_1,\ldots,t_n])$,
$i=1,2$, for all natural numbers $n$.
\end{theorem}

Theorem \ref{BHS} is true for all indices
$i=0,1,2,\ldots$\footnote{and for noncommutative regular rings as
well.} The case $i=0$ is due to Grothendieck, the case $i=1$ is
due to Bass-Heller-Swan \cite{BaHS}, and the general case $i\ge2$
is due to Quillen \cite{Q1}.

\begin{theorem}[\cite{Su}]\label{Suslin}
Let $R$ be a noetherian ring with $\dim R<\infty$ and $n$ be a
nonnegative integer. Then the natural homomorphisms
$$
K_{1,r}(R[t_1,\ldots,t_n])\to K_1(R[t_1,\ldots,t_n])
$$
are surjective for $r\ge\max(2,\dim R+1)$ and bijective for
$r\ge\max(3,\dim R+2)$.
\end{theorem}

Theorems \ref{BHS} and \ref{Suslin} have the following

\begin{corollary}\label{SLN}
Let $\k$ be a field and $n$ be a natural number. Then
$$
\SL_r(\k[t_1,\ldots,t_n])=\E_r(\k[t_1,\ldots,t_n]),\qquad r\ge3.
$$
\end{corollary}

Suslin proves this equality in \cite{Su} directly, without
invoking the Bass-Heller-Swan isomorphism. This is done by
developing a $K_1$-analog of \emph{Quillen's local-global
patching} and \emph{Horrocks' monic inversion technique}, the two
crucial ingredients in Quillen's proof of Serre's conjecture on
projective modules \cite{Q2}. It is exactly Suslin's proof of
Corollary \ref{SLN} what is used in the algorithm, developed in
\cite{PW}:

\begin{theorem}[\cite{PW}]\label{PW}
Let $\k$ be a field and $n$ be a natural number. There is an
algorithm which for any matrix $A\in\SL_r(\k[t_1,\ldots,t_n])$
finds a factorization of the form:
$$
A=\prod_ke_{p_kq_k}(\lambda_k),\qquad
\lambda_k\in\k[t_1,\ldots,t_n].
$$
\end{theorem}

\begin{remark}\label{Cohn}
The inequality $r\ge3$ is sharp as shown by the following example
of Cohn \cite{Coh}. For any field $\k$ we have
$$
A=\begin{pmatrix}
1+t_1t_2&-t_1^2\\
\\
t_2^2&1-t_1t_2
\end{pmatrix}\in\SL_2(\k[t_1,t_2])\setminus\E_2(\k[t_1,t_2]).
$$
By Corollary \ref{SLN}, $A$ becomes an elementary matrix already
in $\SL_3(\k[t_1,t_2])$. However, if we consider the monomial ring
$\k[t_1^2,t_1t_2,t_2^2]$ over which $A$ is defined, then the
matrix $A$ represents a non-zero element in
$K_1(\k[t_1^2,t_1t_2,t_2^2])$, \cite[Example 8.2]{G3}. Therefore,
$A$ does not become an elementary matrix in any of the groups
$\SL_r(\k[t_1^2,t_1t_2,t_2^2])$, no matter how large $r$ is. This
explains the relevance of Frobenius actions (that is, the
homomorphisms $c_*$) in the nilpotence conjecture.
\end{remark}

\begin{remark}\label{whynilp} For a field $\k$ one can sandwich the
2-dimensional polynomial rings between two copies of
$\k[t_1^2,t_1t_2,t_2^2]$ as follows
$$
\k[t_1^2,t_1t_2,t_2^2]\subset \k[t_1,t_2]\subset
\k[t_1,t_1^{1/2}t_2^{1/2},t_2]\cong\k[t_1^2,t_1t_2,t_2^2].
$$
This observation and Corollary \ref{SLN} show that
$(2_*)(A)\in\E_r(\k[t_1^2,t_1t_2,t_2^2])$ for all $r\ge3$. An
elaborated version of this argument, in combination with an
excision technique, implies Theorem \ref{MAIN}(b,c) in the special
case when $M$ is a \emph{simplicial monoid,} which means
$M\subset\ZZ^n$ is a finitely generated additive submonoid and the
cone in $\RR^n$ spanned by $M$ is simplicial; see Corollaries
\ref{simplicialalg} and \ref{simplicialK1} below. However, the
existence of such a sandwiched polynomial ring
$$
\k[M]\subset \k[\ZZ_+^n]\subset \k[M^{1/c}],\qquad
M^{1/c}=\{m^{1/c}\ |\ m\in
M\}\subset\ZZ\left[\frac1c\right]\otimes\gp(M)
$$
implies that $M$ is simplicial. This partly explains why the
general case of the nilpotence conjecture is essentially more
difficult than the simplicial case.
\end{remark}

Tulenbaev's result below, proved in \cite{T}, is a $K_2$-analog of
Suslin's work \cite{Su}.

\begin{theorem}\label{Tulenbaev}
Let $R$ be a noetherian ring of finite Krull dimension $\dim R$
and $n$ a natural number. Then the natural homomorphisms
$$
K_{2,r}(R[t_1,\ldots,t_n])\to K_2(R[t_1,\ldots,t_n])
$$
are surjective for $r\ge\max(4,\dim R+2)$ and bijective for
$r\ge\max(5,\dim R+3)$.
\end{theorem}

Earlier van der Kallen had shown that $K_{2,r}(R)=K_2(R)$ for
$r\ge\dim R+3$ \cite{K1}. Correspondingly, we will always write
$K_2(R)$ instead of $K_{2,r}(R)$ when $r$ is as in Theorem
\ref{Tulenbaev}.

\section{Monoids and cones}\label{MonCon}

Here is a quick summary of the generalities on cones and monoids.
For more detailed account the interested reader is referred to
\cite[Chapters 1,\ 2]{BrG}.

\subsection{Polytopes and cones}\label{polytopesandcones}
A \emph{polytope} $P\subset\RR^n$ means the convex hull of
finitely many points in $\RR^n$. This is the same as a compact
intersection of finitely many affine half-spaces in $\RR^n$. For a
polytope $P\subset\RR^n$ its relative interior will be denoted by
$\inte(P)$. A polytope $P\subset\RR^n$ is called \emph{rational}
if it is spanned by rational points. A polytope $P$ is rational if
and only if it is a compact intersection of finitely many affine
half-spaces whose boundaries are rational affine hyperplanes. A
polytope is a \emph{simplex} if it is the convex hull of an
\emph{affinely independent system} of points.

\smallskip The set of nonnegative reals is denoted by $\RR_+$.
For a subset $X\subset\RR^n$ we will use the notation
$\RR_+X=\{\sum_i a_ix_i\ |\ a_i\in\RR_+,\ x_i\in X\}$.

\smallskip A \emph{cone} $C\subset\RR^n$ means a
 subset
of the form $\RR_+X\subset\RR^n$ where $X$ is finite. This is the
same as the intersection of a finite family of halfspaces in
$\RR^n$ whose boundary hyperplanes are linear subspaces of
$\RR^n$. When $X\subset\QQ^n$ (equivalently, the mentioned
halfspaces have rational boundary hyperplanes) the cone is called
\emph{rational}. A cone is \emph{pointed} if it contains no pair
of opposite nonzero vectors. A cone $C\subset\RR^n$ can be
embedded (via a linear map) in $\RR^{\dim C}$. If $C$ is rational
then such an embedding can be chosen to be rational. Further, a
cone $C\subset\RR^n$ is pointed if and only if it can be embedded
in the positive orthant $\RR_+^{\dim C}$.

\smallskip \emph{All our cones will be
assumed to be pointed.}

\smallskip Let $C\subset\RR^n$ be a cone and $\mathcal
H^+\subset\RR^n$ be a half-space, defined by an inequality
$\xi_1X_1+\cdots+\xi_nX_n\ge0$, such that $C\subset\mathcal H^+$.
Let $\mathcal H$ be the boundary hyperplane
$\xi_1X_1+\cdots+\xi_nX_n=0$. Then the intersection $C\cap\mathcal
H$ is called a \emph{face of $C$}. The origin $0$ and the cone $C$
itself are the smallest and the biggest faces of $C$. A
\emph{facet} of a cone $C\subset\RR^n$ is a maximal proper face,
which is the same as a codimension 1 face. The boundary $\partial
C$ is defined as the union of all proper faces of $C$, and
\emph{the relative interior $\inte(C)$} is defined by
$\inte(C)=C\setminus\partial C$.

\smallskip A \emph{$d$-cone} means a $d$-dimensional cone.

\smallskip An \emph{open cone} in $\RR^n$ of dimension $d$ is by definition
the union of the relative interiors of $d$-cones, forming a nested
system of cones, plus the origin $0$.

\smallskip An \emph{affine cone} means a parallel translate of a
cone.

\smallskip For a rational $d$-cone
$C\subset\RR^n$, $d>0$, there always exists a rational affine
$(n-1)$-dimensional subspace $\G\subset\RR^n\setminus\{0\}$ such
that $C=\RR_+(C\cap\G)$ or, equivalently, $C\cap\G$ is a rational
$(d-1)$-polytope. For such a pair $C$ and $\G$ we write
$\Phi(C)=C\cap\mathcal G$. Further, for a real number $\epsilon>0$
we will use the notation $C(\epsilon)=\RR_+\Phi(C)(\epsilon)$
where $\Phi(C)(\epsilon)$ is the $\epsilon$-neighborhood of
$\Phi(C)$ in $\G$. Thus $C(\epsilon)\subset\RR^n$ is an
$n$-dimensional open cone.

\smallskip  A cone is called \emph{simplicial} if it
is spanned by a system linearly independent vectors, or
equivalently, the polytope $\Phi(C)$ is a simplex.

\subsection{Monoids}\label{Monoids} A \emph{monoid} will always mean a commutative,
cancellative, torsion free monoid. Equivalently, our monoids are
additive submonoids of rational vector spaces.

\smallskip Our blanket assumption on the notation of monoid operation is
that when a monoid is considered inside its monoid ring we use
multiplicative notation. Otherwise we use additive notation.

\smallskip For a monoid $M$ its group of differences will be denoted
by $\gp(M)$. We put $\rank M=\rank\gp(M)$. If a monoid is finitely
generated then it is called \emph{affine}. Thus an affine monoid
is, up to isomorphism, a finitely generated additive submonoid of
$\ZZ^n$. Moreover, whenever appropriate we can without loss of
generality assume that $\gp(M)=\ZZ^n$.

\smallskip A monoid is called \emph{positive} if its group
of invertible elements is trivial. For an affine positive monoid
$M\subset\ZZ^n$ the subset $ \RR_+M\subset\RR^n$ is a rational
cone. A monoid $M$ is called \emph{simplicial} if it is positive,
affine and the cone $\RR_+M$ is simplicial.

\smallskip For an affine positive monoid $M\subset\ZZ^n$, $\rank M>0$, and an
affine hyperplane $\mathcal G\subset\RR^n$ such that
$\RR_+M=\RR_+\left(\RR_+M\cap\mathcal G\right)$, we will use the
notation $\Phi(M)$ for $\RR_+M\cap\mathcal G$. For a convex subset
$W\subset\Phi(M)$ we introduce the submonoid
$$
M|W=M\cap\RR_+W\subset M.
$$
If $W$ consists of a single point $p$ then we write $M|p$ instead
of $M|\{p\}$.

\smallskip For $M$ and $\mathcal G$ as above we will also use the notation
$M_*=M\cap\RR_+\inte(\Phi(M))$\footnote{Here we follow the
convention that the interior of a point is the point itself. In
particular, $M=M_*$ when $\rank M=1$.} and $M|F=M\cap F\subset M$
for $F\subset\RR_+M$ a face. Thus
$M_*=\big(M\cap\inte(\RR_+M)\big)\cup\{0\}$. More generally, if
$N\subset M$ is any (not necessarily affine) submonoid then we put
$\Phi(N)=\mathcal G\cap \RR_+N$ and
\begin{align*}
N_*=N\cap\RR_+\{x\ |\ x\ \text{is in the relative interior of}\
\Phi(N)\}.
\end{align*}

\smallskip For an affine positive monoid $M\subset\ZZ^n$ and a convex subset
$W\subset\Phi(M)$ (w.r.t. to an appropriately fixed hyperplane
$\G\subset\RR^n$ as above) it is easily shown that
\begin{equation}\label{gpequal}
\dim W=\rank M-1\ \implies\ \gp(M)=\gp(M|W).
\end{equation}
(See, for instance, \cite[Corollary 2.25]{BrG}.) In particular,
\begin{equation}\label{gpm*}
\gp(M)=\gp(M_*).
\end{equation}

\smallskip Let $M\subset\ZZ^n$ be an affine positive monoid,
$F\subset\RR_+M$ a face, and $R$ a ring. Then we have the
$R$-algebra retraction:
\begin{align*}
\pi_F:R[M]\to R[M|F],\quad \pi(m)=
\begin{cases}
m\ \text{if}\ m\in M|F,\\
0\ \text{if}\ m\in M\setminus(M|F).
\end{cases}
\end{align*}

\smallskip A monoid $M$ is called \emph{normal} if
$kx\in M$ implies $x\in M$ for any $x\in\gp(M)$ and any $k\in\NN$.
Any affine positive normal monoid of rank $n$ is up to isomorphism
of the form $C\cap\ZZ^n$ where $C\subset\RR^n$ is a positive
rational $n$-cone. Conversely, any such an intersection
$C\cap\ZZ^n$ is always an affine positive normal monoid. The
finite generation part of the latter claim is classically known as
\emph{Gordan's lemma} (\cite[Lemma 2.7]{BrG}).

\smallskip For any monoid $M$ there is the smallest submonoid of
$\gp(M)$ -- the \emph{normalization of $M$} -- which is normal and
contains $M$:
$$
\bar M=\{x\in\gp(M)\ |\ kx\in M\ \text{for some natural number}\
k\}.
$$

\smallskip For an affine normal positive monoid $M\subset\ZZ^n$
and a convex subset $W\subset\mathcal G$, where $\mathcal
G\subset\RR^n$ is a hyperplane cross-secting $\RR_+M$, we
introduce the monoid:
$$
M|W=\gp(M)\cap\RR_+W.
$$
When $W\subset\Phi(M)$ this notation is compatible with the one
introduced above for not necessarily normal monoids.

\smallskip A monoid $M$ is called \emph{seminormal} if the following
implications holds:
$$
x\in\gp(M),\ 2x\in M,\ 3x\in M\ \implies x\in M.
$$

\begin{lemma}\label{smn}
Let $M\subset\ZZ^n$ be an affine positive monoid. Then $M$ is
seminormal if and only if the monoid $(M|F)_*$ is normal for any
face $F\subset\RR_+M$. Moreover, if $M$ is seminormal then
$M_*=\bar M_*$.
\end{lemma}

The first part is proved in \cite{G1} (for not necessarily affine
monoids), see also \cite[Proposition 2.37]{BrG}. The second part
follows from the equality (\ref{gpm*}).

\subsection{Divisible monoids}\label{seminormal} For a natural number $c$ and a monoid
$M$ we say that $M$ is \emph{$c$-divisible} if for any element
$z\in M$ the equation $cx=z$ is solvable for $x$ inside $M$. Since
our monoids are cancellative and torsion free, such a solution is
unique.

\smallskip For a monoid $M$ and a natural number $c$ the submonoid
of $\ZZ\left[\frac1c\right]\otimes\gp(M)$, generated by
$\frac1c\otimes x$, $x\in M$, will be denoted by $M/c$.

\smallskip For a natural number $c\ge2$ the \emph{$c$-divisible
hull} of $M$ is defined as the filtered union
$$
M/{c^{\infty}}=\bigcup_{j=1}^\infty M/c^j\subset\QQ\otimes\gp(M).
$$

It is easily checked that for a natural number $c\ge2$ all $c$
divisible monoids $L$ are seminormal:
\begin{align*}
2x,3x\in L\ \Longrightarrow\ cx\in L\ \Longrightarrow
x=\frac1c\cdot(cx)\in L.
\end{align*}

By Lemma \ref{smn} the submonoid $M_*/{c^{\infty}}\subset
M/{c^{\infty}}$ is a normal monoid for any positive affine monoid
$M$. It easily follows that for any affine positive monoid $M$ we
have:
\begin{equation}\label{M*}
(M_*)^{c^{-\infty}}=(\bar M_*)^{c^{-\infty}}.
\end{equation}

\smallskip When $M$ is simplicial much more is true:

\begin{proposition}\label{apprA}
Let $M$ be an affine simplicial monoid. Then for any finite subset
$S\subset M_*/{c^{\infty}}$ one can effectively find a free
submonoid $L\subset M_*/{c^{\infty}}$ such that $S\subset L$. In
particular, $M_*/{c^{\infty}}$ is a filtered union of free
monoids.
\end{proposition}

Without effective nature of the claim this is Theorem A in
\cite{G2}. However, what is proved in \cite{G2} is literally what
is stated above.

\

Next we derive a structural result on $c$-divisible monoids that
will be used in Section \ref{lastsection}. Let $M\subset\ZZ^n$ be
an affine positive monoid and let $h:\ZZ^n\to\ZZ$ be a surjective
group homomorphism. Then $M$ carries the \emph{graded structure:}
\begin{align*}
M=\cdots\cup M_{-1}\cup M_0\cup M_1\cup\cdots,\qquad M_i=M\cap
h^{-1}(i).
\end{align*}
(`Graded' here means $M_i+M_j\subset M_{i+j}$ and $M_i\cap
M_j=\emptyset$ whenever $i\not=j$.) For an element $m\in M_i$ we
will write $\deg(m)=i$.

\smallskip For simplicity of notation we let the same $h$ denote the
$\RR$-linear extension $\RR\otimes h:\RR^n\to\RR$.

\begin{lemma}\label{decomposition}
Let $M\subset\ZZ^n$ be an affine positive monoid with
$\gp(M)=\ZZ^n$. Let $m\in M_*$ with $\deg(m)=d\not=0$. Then one
can effectively find a decomposition of the form:
\begin{align*}
m=\sum_{i=1}^{|d|} m_i,\quad m_i\in
\begin{cases}
(M_*)^{c^{-\infty}}\cap h^{-1}(1) &\text{if}\ d>0,\\
(M_*)^{c^{-\infty}}\cap h^{-1}(-1) &\text{if}\ d<0.
\end{cases}
\end{align*}
\end{lemma}

\begin{proof}
We consider the case $d>0$ and the other case is symmetric.

Consider the broken line $\mathfrak a=[0a_1a_2\ldots a_{d-1}m]$ in
$\RR^n$, obtained by subdividing the segment $[0,m]\subset\RR^n$
into $d$ equal parts. This broken line can be though of as the
decomposition inside $\QQ^n$:
$$
m=d^{-1}m+\cdots+d^{-1}m.
$$
We want to find (effectively!) a broken line in $\RR^n$
$$
\mathfrak m=[m_0m_1m_2\ldots m_d],\quad m_0=0,\quad m_d=m,
$$
satisfying the condition $m_i-m_{i-1}\in(M_*)^{c^{-\infty}}\cap
h^{-1}(1)$ for $i=1,\ldots,d-1$. Since $\RR_+M_*$ is an open cone,
any broken line $\mathfrak b=[m_0m_1m_2\ldots m_{d-1}m_d]$ that is
obtained from $\mathfrak a$ by an arbitrary sufficiently small
perturbation of the vertices $a_1,\ldots,a_{d-1}$ will satisfy the
condition $m_i-m_{i-1}\in\RR_+M_*$ for $i=1,\ldots,d-1$.
Therefore, it is enough to show that for every index
$i\in\{1,\ldots,d-1\}$ the affine real hyperplane
$h^{-1}(i)\subset\RR^n$ contains elements of $(M_*)^{c^{-\infty}}$
arbitrarily close to $a_i$. In view of the equalities
(\ref{gpequal}) and (\ref{M*}),  it is enough to show that for
every index $i\in\{1,\ldots,d-1\}$ the affine real hyperplane
$h^{-1}(i)\subset\RR^n$ contains elements of
$\gp(M)^{c^{-\infty}}$ arbitrarily close to $a_i$. This will be
done by showing that for every $i\in\{1,\ldots,d-1\}$ the set
$\gp(M)^{c^{-\infty}}\cap h^{-1}(i)$ is dense in the affine
hyperplane $h^{-1}(i)\subset\RR^n$.

The conditions $\gp(M)=\ZZ^n$ and $h(\ZZ^n)=\ZZ$ imply that the
sets
\begin{align*}
h^{-1}(i)\cap\gp(M),\quad i=1,\ldots,d-1,
\end{align*}
are cosets of $\Ker(h)\cap\ZZ^n$ in $\ZZ^n$. In particular,
\begin{align*}
\gp(M)^{c^{-\infty}}\cap
h^{-1}(i)\simeq(\ZZ[1/c])^n\cap\Ker(h),\quad i=1,\ldots,k-1,
\end{align*}
where $\ZZ[1/c]$ refers to the localization of the ring of
integers $\ZZ$ at $c$ and $\simeq$ refers to the isometry
equivalence w.r.t. the Euclidean metric. But
$(\ZZ[1/c])^n\cap\Ker(h)$ is a $c$-divisible $\rank(n-1)$ subgroup
of $\Ker(h)\cong\RR^{n-1}$. In particular, it is a dense subset of
$\Ker(h)$.

\smallskip The algorithmic aspect of Lemma \ref{decomposition}
follows from the fact that we can effectively compute (in terms of
generators) the group $\ZZ^n\cap\Ker(h)$, its appropriate cosets
in $\ZZ^n$, and find an element of
$\gp(M)^{c^{-\infty}}\cap\Ker(h)$ in any explicitly given
neighborhood in $\Ker(h)$.
\end{proof}

\smallskip The multiplicative counterpart of the notation
$M/c^j$ and $M/{c^{\infty}}$, to be used in monoid rings, is
$M^{c^{-j}}$ and $M^{c^{-\infty}}$.

\smallskip The relevance of $c$-divisible monoids is explained by
the following equivalent reformulation of Theorem \ref{MAIN}:

\begin{theorem}\label{cdivrev}
Let $M$, $c$, $R$ and $\k$ be as in Theorem \ref{MAIN}. Then
\begin{itemize} \item[(a)]
$K_2(R)=K_{2,r}(R[M^{c^{-\infty}}])$ for $r\ge\max(5,\dim R+3)$.
\item[(b)]
$\GL_r(R[M^{c^{-\infty}}])=\E_r(R[M^{c^{-\infty}}])\GL_r(R)$ for
$r\ge\max(3,\dim R+2)$. \item[(c)] There is an algorithm which for
any matrix $A=\SL_r(\k[M])$, $r\ge3$, finds an integer number
$j_A\ge0$ and a factorization of the form:
$$
A=\prod_k e_{p_kq_k}(\lambda_k),\qquad
e_{p_kq_k}(\lambda_k)\in\E_r(\k[M^{c^{-j_A}}]).
$$
\end{itemize}
\end{theorem}
In the subsequent sections we will freely use the equivalence
between the two formulations.

\begin{remark}\label{WHY?}
Essentially, $c$-divisible monoids enter our argument through
Proposition \ref{apprA} (and a variation of it -- Lemma
\ref{filtun}) and Lemma \ref{decomposition}, used correspondingly
in Sections \ref{redint} and \ref{Mushk}. They also partially
explain why in this paper we mainly work with open cones. In
\cite{G5} the importance of $c$-divisible monoids is related to
the excision results in \cite{SuW} and that of open cones -- to
Karoubi squares of certain type.
\end{remark}

\section{Reduction to interior monoids}\label{redint}

\begin{proposition}\label{efectmiln}
Let $M\subset\ZZ^n$ be an affine positive monoid. Assume Theorem
\ref{MAIN} is valid for the submonoids of the form $(M|F)_*\subset
M$ where $F\subset\RR_+M$ is a facet or $F=\RR_+M$. Then the
theorem is valid also for $M$.
\end{proposition}

For a matrix $A\in\GL_r(R[M])$ the elements $m\in M$ that show up
in the canonical $R$-linear expansion of its entries will be
called the \emph{support monomials of $A$}.

\smallskip A monoid is a filtered union of its affine submonoids. Moreover,
one can find effectively such a filtered union representation for
any explicitly given monoid. Therefore, by the equality (\ref{M*})
in Section \ref{seminormal}, Proposition \ref{efectmiln} and the
equivalent reformulation of Theorem \ref{MAIN} in Theorem
\ref{cdivrev} we get

\begin{corollary}\label{internormred}
For Theorem \ref{MAIN} it is enough to show that for any affine
positive normal monoid $M$ we correspondingly have:
\begin{itemize}
\item[(a)] $K_2(R)=K_{2,r}(R[(M_*)^{c^{-\infty}}])$, \item[(b)]
$\GL_r(R[(M_*)^{c^{-\infty}}])=\E_r(R[(M_*)^{c^{-\infty}}])\GL_r(R)$,
\item[(c)] There is an algorithm that for any matrix
$A\in\SL_r(\k[M_*])$ finds an integer number $j_A\ge0$ and a
factorization of the form:
\begin{align*}
A=\prod_ke_{p_kq_k}(\lambda_k),\qquad
e_{p_kq_k}(\lambda_k)\in\E_r(\k[(M_*)^{c^{-j_A}}]).
\end{align*}
\end{itemize}
(Here $r$ is as in the corresponding part of Theorem \ref{MAIN}.)
\end{corollary}

\smallskip In the next three subsections we prove
Proposition \ref{efectmiln}, considering the three parts of
Theorem \ref{MAIN} separately and in the reversed order. The case
of Milnor groups requires substantially more work.

\subsection{The case of Theorem
\ref{MAIN}(c)}\label{inter-c} Let $F\subset \RR_+M$ be a facet and
$A\in\SL_r(\k[M])$. Consider the matrix
$A|F=\pi_F(A)\in\SL_r(\k[M|F])$. Obviously, $A|F$ is effectively
computable from $A$: its support monomials are those of $A$ that
belong to $F$. By the assumption, $A|F$ can be effectively
factored into standard elementary matrices over
$\k[(M|F)^{c^{-j_F}}]$ for some explicitly computable $j_F\in\NN$.
Therefore, it is enough to prove Theorem \ref{MAIN}(c) for the
matrix $A_F=(A|F)^{-1}A\in\SL_r(\k[M])$. Observe that no support
monomial of $A_F$ belongs to $M|F$.

Now let $G\subset\RR_+M$ be another facet. Again by the assumption
the matrix $A_F|G=\pi_G(A_F)\in\SL_r(\k[M|G])$ can be
algorithmically factored into elementary matrices over the ring
$\k[(M|G)^{c^{-j_G}}]$ for some explicitly computable $j_G\in\NN$.
It is enough to prove Theorem \ref{MAIN}(c) for the matrix
$A_{F,G}=(A_F|G)^{-1}A_G\in\SL_r(\k[M])$.

The crucial observation at this point is that no support monomial
of the matrix $A_{F,G}$ belongs to $(M|F)\cup(M|G)$.

Continuing the process until all facets of the cone $\RR_+M$ are
considered, we arrive at a matrix
$$
A_{F,G,\ldots,H}\in\SL_r(\k[M_*])
$$
where $\{F,G,\ldots,H\}$ is the set of facets of $\RR_+M$. By the
assumptions in the proposition, one can find $j_M\in\NN$ and a
factorization of $A_{F,G,\ldots,H}$ into standard elementary
matrices from $\E_r(\k[(M_*)^{c^{-j_M}}])$.

It is then clear that the desired explicit factorization of $A$
can be found over the ring
$\k[M_*^{c^{-j_F-j_G-\cdots-j_H-j_M}}]$. \qed

\smallskip In view of Theorem \ref{PW} and Proposition \ref{apprA}
the induction on $\rank M$ yields

\begin{corollary}\label{simplicialalg}
Theorem \ref{MAIN}(c) is true for simplicial monoids.
\end{corollary}

\subsection{The case of Theorem \ref{MAIN}(b)} Essentially the
same argument as above goes through. In more detail, consider a
matrix $A\in\GL_r(R[M])$. For a facet $F\subset\RR_+M$ we have the
matrix $A|F=\pi_F(A)\in\GL_r(R[M|F])$. By the assumptions in the
proposition, there exists
$E_F\in\E_r(R[(M|F)^{c^{-\infty}}])\subset\E_r(R[M^{c^{-\infty}}])$
such that $E_F\cdot(A|F)\in\GL_r(R)$. In particular,
$E_F\in\GL_r(R[M|F])$ and no support monomial of the matrix $E_FA$
belongs to $M|F$.

It is enough to show that
$E_FA\in\E_r(R[M^{c^{-\infty}}])\GL_r(R)$.

Consider another facet $G\subset\RR_+M$. Again by the induction
hypothesis there exists
$E_G\in\E_r(R[(M|G)^{c^{-\infty}}])\subset\E_r(R[M^{c^{-\infty}}])$
such that $E_G\cdot((E_FA)|G)\in\GL_r(R)$. In this situation
$E_G\in\GL_r(R[M|G])$ and no support monomial of $E_GE_FA$ belongs
to $M|G$. We claim that no support monomial of $E_GE_FA$ belongs
to $M|F$ too. In fact, we have
$\pi_F(E_GE_FA)=\pi_F(E_G)\pi_F(E_FA)\subset\GL_r(R[M|G])\GL_r(R)$.
In particular, if there were a support monomial of $E_GE_FA$ in
$M|F$ then it would also belong to $M|G$. But such does not exist.

Continuing the process with the remaining facets we find a system
of elementary matrices
\begin{align*}
E_F\in\E_r(R[(M|F)^{c^{-\infty}}]),\
E_G\in\E_r(R[(M|G&)^{c^{-\infty}}]),
\ldots,E_H\in\E_r(R[(M|H)^{c^{-\infty}}]),\\
&F,G,\ldots,H\subset\RR_+M\ \text{are the facets},
\end{align*}
such that $ E_H\cdots E_GE_FA\in\GL_r(R[M_*])$. But over $R[M_*]$
we are done by the assumptions in the proposition. \qed

\smallskip In view of Theorems \ref{BHS}, \ref{Suslin} and Proposition
\ref{apprA} the induction on $\rank M$ yields

\begin{corollary}\label{simplicialK1}
Theorem \ref{MAIN}(b) is true for simplicial monoids.
\end{corollary}

\subsection{The case of Theorem \ref{MAIN}(a)} This is not as straightforward as the
previous cases.

\begin{lemma}\label{filtun}
Let $c\ge2$ be a natural number and $M_1,M_2$ be $c$-divisible
monoids of rank 1 without nontrivial units. Then the submonoid
$$
N=M_1\times M_2\setminus\{(a,0)\ |\ a\in M_1,\ a\not=0\}\subset
M_1\times M_2
$$
is a filtered union of rank 2 free monoids.
\end{lemma}

\begin{proof}
There are inductive systems of indices $I$ and $J$ and elements
$a_i\in M_1$, $i\in I$, and $b_j\in M_2$, $j\in J$, such that:

\smallskip$\centerdot$ $M_1=\bigcup_I A_i$ and $M_2=\bigcup_JB_j$,

\smallskip$\centerdot$ $A_i=\ZZ_+a_i$ and if $i_1<i_2$ then
$s_{i_2,i_1}a_{i_2}=a_{i_1}$ for some natural number
$s_{i_2,i_1}\ge2$,

\smallskip$\centerdot$ $B_j=\ZZ_+b_j$ and if $j_1<j_2$ then
$t_{j_2,j_1}b_{j_2}=b_{j_1}$ for some natural number
$t_{j_2,j_1}\ge2$.

\

For any pair $(i,j)\in I\times J$ consider the monoid
$$
N_{ij}=\ZZ_+(a_i,b_j)+\ZZ_+(0,b_j)\cong\ZZ_+^2
$$
For any indices $i\in I$ and $j_1,j_2\in J$ with $j_1\le j_2$ we
have
\begin{align*}
(a_i,b_{j_1})=(a_i,b_{j_2})+(t_{j_2,j_1}-1)(1,b_{j_2})\in
N_{ij_2}.
\end{align*}
Therefore, $N_{ij_1}\subset N_{ij_2}$. In particular, the monoids
$$
N_i=A_i\times M_2\setminus\{(a,0)\ |\ a\in A_i,\ a\not=0\}
$$
are filtered unions of the monoids $N_{ij}$, $j\in J$. But $N$ is
a filtered union of the monoids $N_i$.
\end{proof}

In the next lemma we use the following notation: for a
homomorphism of rings $\Lambda_1\to\Lambda_2$ and a natural number
$r$ we let $\st^*_r(\Lambda_1)$ denote the image of the map
$\st_r(\Lambda_1)\to\st_r(\Lambda_2)$.

\begin{lemma}\label{segment}
Let $R$ be a regular ring of finite Krull dimension $d$ and
$r\ge\max(5,d+3)$. Assume $c$, $M_1$, $M_2$ and $N\subset
M_1\times M_2$ are as in Lemma \ref{filtun}. Then arbitrary
element $w\in\st_r(R[M_1\times M_2])$ admits a presentation of the
form
\begin{align*}
w=uv,\quad u\in\st^*_r(R[M_1]),\quad v\in\st^*_r(R[N]).
\end{align*}
\end{lemma}
\noindent(Here the maps from $\st_r(R[M_1])$ and $\st_r(R[N])$ to
$\st_r(R[M_1\times M_2])$ are the ones induced by the identity
ring embeddings $R[M_1]\to R[M_1\times M_2]$ and $R[N]\to
R[M_1\times M_2]$.)

\begin{proof}
Consider the commutative square of $R$-algebra homomorphisms whose
horizontal arrows are identity embeddings:
$$
\xymatrix{
R[N]\ar[r]\ar[d]_{\theta|_{R[N]}}&R[M_1\times M_2]\ar[d]^\theta\\
R\ar[r]&R[M_1] },\quad\theta|_{M_1}=\1_{M_1},\quad
\theta(M_2\setminus\{1\})=0.
$$

Because $\st_r(R)\to\st_r(R[M_1])\to\st_r(R[M_1\times M_2])$ are
(split) injective homomorphisms, we can identify $\st_r(R)$ and
$\st_r(R[M_1])$ with the subgroups
$\st^*_r(R)\subset\st^*_r(R[M_1])\subset\st_r(R[M_1\times M_2])$.

Let $\tau:\st^*_r(R[N])\to\st_r(R)$ be the homomorphism induced by
the augmentation $R[M_1\times M_2]\to R$, $M_1\times
M_2\setminus\{(1,1)\}\to0$.

\smallskip First we show the following inclusion
\begin{equation}\label{kerst}
\Ker(\tau)\st_r(R[M_1])\subset\st_r(R[M_1])\Ker(\tau).
\end{equation}

\smallskip Assume $u_1\in\st_r(R[M_1])$ and $v_1\in\Ker(\tau)$. We
want to prove that $v_1u_1\in\st_r(R[M_1])\Ker(\tau)$.

\smallskip Let $v'=u_1^{-1}v_1u_1\in\st_r(R[M_1\times M_2])$ and
$e_1$ and $e'$ be the images of $v_1$ and $v'$ in
$\E_r(R[M_1\times M_2])$. From the commutative square
$$
\xymatrix{ \st^*_r(R[N])\ar[r]^{\tau}\ar[d]&\st_r(R)\ar[d]\\
\E_r(R[M_1\times M_2])\ \ar[r]^{\E_r(\theta)}&\E_r(R[M_1])
 }
$$
we see that $e_1\in\Ker(\E_r(\theta))$. Then
$e'\in\Ker(\E_r(\theta))$ as well. In particular,
$e'\in\SL_r(R[N])$.

By Lemma \ref{filtun} $N$ is a filtered union of rank 2 monoids.
Therefore, by Theorems \ref{BHS} and \ref{Suslin} we have
$\GL_r(R[N])=\GL_r(R)\E_r(R[N])$ and so
$$
e'\in\Ker(\E_r(\theta))\cap\GL_r(R)\E_r(R[N])\subset\E_r(R[N]).
$$

Let $v''\in\st^*_r(R[N])$ be a preimage of $e'$. There exists
$z\in K_{2,r}(R[M_1\times M_2])$ such that $v'=zv''$. The monoid
$M_1\times M_2$ is clearly a filtered union of rank 2 free monoids
and so $K_{2,r}(R[M_1\times M_2])=K_2(R)$ by Theorems \ref{BHS}
and \ref{Tulenbaev}. Hence the desired representation
\begin{align*}
v_1u_1=u_2v_2,\quad u_2=u_1z\tau(v'')\in\st_r(R[M_1]),\quad
v_2=\tau(v'')^{-1}v''\in\Ker(\tau).
\end{align*}

\smallskip Finally, Lemma \ref{segment} follows from
(\ref{kerst}) because any generator $x_{ij}(\lambda)$ of the group
$\st_r(R[M_1\times M_2])$ has a representation of the form:
\begin{align*}
x_{ij}(\lambda)=x_{ij}(\theta(\lambda))x_{ij}(\lambda-\theta(\lambda)),\
x_{ij}(\theta(\lambda))\in\st_r(R[M_1]),\
x_{ij}(\lambda-\theta(\lambda))\in\Ker(\tau).
\end{align*}
\end{proof}

From now on we assume that $c$, $R$, $r$ and $M$ are as in Theorem
\ref{MAIN}(a).

\smallskip Fix a facet $F\subset\Phi(M)$. By the induction
hypothesis we have
\begin{equation}\label{facetdone}
K_{2,r}(R[(M|F)^{c^{-\infty}}])=K_2(R).
\end{equation}

\smallskip Any element $z\in K_{2,r}(R[M])$ has a representation
of the form $z=\prod_k v_k$ where:

\smallskip$\centerdot$ $v_{k_1}\in\st^*_r(R[(M|p_{k_1})^{c^{-\infty}}]),\ldots,v_{k_s}\in
\st^*_r(R[(M|p_{k_s})^{c^{-\infty}}])$,

\smallskip$\centerdot$ $v_k\in\st^*_r(R[(M|q_k)^{c^{-\infty}}])$, $k\notin\{k_1,\ldots,k_s\}$,

\smallskip$\centerdot$ $k_1<\ldots<k_s,\quad p_{k_1},\ldots,p_{k_s}\in
F$, $q_k\in\Phi(M)\setminus F$, $k\notin\{k_1,\ldots,k_s\}$.

\smallskip\noindent(For instance, any representation of the form $\prod_k
x_{i_kj_k}(\mu_k)$ where $\mu_k\in R M^{c^{-\infty}}$ is of this
form.)

\smallskip When $\{k_1,\ldots,k_s\}\not=\emptyset$ we say that
$z$ has a representation of \emph{$(k_1,\ldots,k_s)$-type}.

\begin{lemma}\label{1stype}
If $z\in K_{2,r}(R[M^{c^{-\infty}}])$ has a representation of
$(1,\ldots,s)$-type then
$$
z\in\Im(K_{2,r}(R[(M|\Phi(M)\setminus F)^{c^{-\infty}}])\to
K_{2,r}(R[M^{c^{-\infty}}])).
$$
\end{lemma}

\begin{proof}
Let $z=\prod_kv_k$ be a representation of $(1,\ldots,s)$-type.
Then, denoting by $e_k\in\E_r(R[M^{c^{-\infty}}])$ the image of
$v_k$, $k=1,\ldots,s$, we have
\begin{align*}
&\prod_{k\le s}e_k=\left(\prod_{k>s}
e_k\right)^{-1}\in\E_r(R[(M|F)^{c^{-\infty}}])
\cap\E_r(R[(M|\Phi(M)\setminus
F)^{c^{-\infty}}])\subset\\
&\SL_r(R)\cap\E_r(R[(M|F)^{c^{-\infty}}])=\E_r(R)
\end{align*}
(The latter equality follows from the fact that $R$ is a retract
of $R[(M|F)^{c^{-\infty}}]$.) In particular, there exists an
element $z_1\in\st_r(R)$ such that $zz_1^{-1}\in
K_{2,r}(R[(M|F)^{c^{-\infty}}])=K_2(R)$ (by (\ref{facetdone})).
Now the lemma follows because
\begin{align*}
z=\left(zz_1^{-1}\right)&\left(z_1\prod_{k>s}
v_k\right)\in\Im(K_{2,r}(R[(M|\Phi(M)\setminus
F)^{c^{-\infty}}])\to K_{2,r}(R[M^{c^{-\infty}}])).
\end{align*}
\end{proof}

\begin{lemma}\label{ktol}If $z\in K_{2,r}(R[M^{c^{-\infty}}])$
has a representation of $(k_1,\ldots,k_s)$-type for some
$(k_1,\ldots,k_s)\not=(1,\ldots,s)$ then $z$ has a representation
of $(l_1,\ldots,l_s)$-type for some
$(l_1,\ldots,l_s)<(k_1,\ldots,k_s)$ w.r.t. the lexicographical
order.
\end{lemma}

\begin{proof} Let $z=\prod_k v_k$ be a representation of
$(k_1,\ldots,k_s)$-type and $i\in\{1,\ldots,s\}$ be the smallest
index with $i< k_i$. Thus
$(k_1,\ldots,k_s)=(1,2,\ldots,i-1,k_i,k_{i+1},\ldots,k_s)$. (We do
not exclude the case when $i=1$.)

In this situation we have
$v_{k_i-1}\in\st^*_r(R[(M|q)^{c^{-\infty}}])$ for some
$q\in\Phi(M)\setminus F$ and
$v_{k_i}\in\st^*_r(R[(M|p)^{c^{-\infty}}])$ for some $p\in F$. By
Lemma \ref{segment} we can write
\begin{align*}
z=\left(\prod_{k<k_i-1}v_k\right)\cdot
(v_{k_i-1}v_{k_i})\cdot\left(\prod_{k\ge
k_i+1}v_k\right)=\left(\prod_{k<k_i-1}v_k\right)\cdot
(uv)\cdot\left(\prod_{k\ge k_i+1}v_k\right)
\end{align*}
for some $u\in\st^*_r(R[(M|p)^{c^{-\infty}}])$ and
$v\in\st^*_r(R[\big(M|[q,p)\big)^{c^{-\infty}}])$. Here $[q,p)$
refers to the corresponding half-open segment in $\Phi(M)$.

There exists a representation of the form $v=\prod_jw_j$ where
$w_j\in\st^*_r(R[(M|t_j)^{c^{-\infty}}])$ for some
$t_j\in\Phi(M)\setminus F$. Then
\begin{align*}
z=\left(\prod_{k<k_i-1}v_k\right)\cdot
u\cdot\left(\prod_jw_j\right)\cdot\left(\prod_{k\ge
k_i+1}v_k\right)
\end{align*}
is a representation of
$(1,2,\ldots,i-1,k_i-1,k'_{i+1},\ldots,k_s')$-type for some
$k'_{i+1}\ge k_{i+1}$, $\ldots$, $k'_s\ge k_s$.
\end{proof}

By Lemmas \ref{1stype} and \ref{ktol} we have

\begin{corollary}\label{surjk2}
The identity embedding $R[(M|(\Phi(M)\setminus
F))^{c^{-\infty}}]\to R[M^{c^{-\infty}}]$ induces a surjective
homomorphism
$$ \iota_F:K_{2,r}(R[(M|(\Phi(M)\setminus
F))^{c^{-\infty}}])\to K_{2,r}(R[M^{c^{-\infty}}]).
$$
\end{corollary}

\smallskip Now we complete the proof of Proposition \ref{efectmiln} as follows.

\smallskip Consider a facet $F\not=G\subset\Phi(M)$. Applying the
same argument as in the proof of Corollary \ref{surjk2} to the
elements of $\Im(\iota_F)$ we arrive to the conclusion that the
natural homomorphism
$$
\iota_{F,G}:K_{2,r}(R[(M|(\Phi(M)\setminus (F\cup
G)))^{c^{-\infty}}])\to K_{2,r}(R[(M|(\Phi(M)\setminus
F))^{c^{-\infty}}])
$$
is also surjective. Then we consider another facet of $\Phi(M)$
etc. Finally we obtain the surjectivity of the composite
homomorphism
\begin{align*}
\iota_{F,G,\ldots,H}:K_{2,r}(R[(M_*)^{c^{-\infty}}])\to
K_{2,r}(R[M^{c^{-\infty}}]).
\end{align*}
But $K_{2,r}(R[(M_*)^{c^{-\infty}}])=K_2(R)$. \qed

\

In view of Theorems \ref{BHS}, \ref{Tulenbaev} and Proposition
\ref{apprA} the induction on $\rank M$ yields
\begin{corollary}\label{simplK2}
Theorems \ref{MAIN}(a) is true for simplicial monoids.
\end{corollary}

\section{Pyramidal descent}
\label{proofmain} In this section we introduce a polyhedral
induction technique in $K$-theory of monoid rings, called
\emph{pyramidal descent}, here adapted to the situation of Theorem
\ref{MAIN}. It was introduced in \cite{G1} and further refined in
\cite{G5}. We in fact need the refinement of the technique as
developed in \cite{G5}, see Remark \ref{whycompl}.

\subsection{Pyramidal extensions of polytopes}\label{pyrext} A polytope $P\subset\RR^n$ is
called a \emph{pyramid} if it is a convex hull of one of its
facets $F\subset P$ and a vertex $v\in P$, not in the affine hull
of $F$. In this situation $F$ is a \emph{base} and $v$ is an
\emph{apex} of $P$, and we write $P=\pyr(v,F)$. For instance, an
arbitrary simplex is a pyramid such that every facet is a base and
every vertex is an apex.

The \emph{complexity} of a $d$-dimensional polytope
$P\subset\RR^n$ is defined as the number $\cc(P)=d-i$, where $i$
is the maximal nonnegative integer satisfying the condition: there
exists a sequence $P_0\subset P_1\subset\cdots\subset P_i=P$ such
that $P_j$ is a pyramid over $P_{j-1}$ for each $1\leq j\leq i$.
Observe that if $P$ is a rational polytope then so are the
polytopes $P_0, P_1,\ldots,P_{i-1}$.

Informally, the complexity of a polytope is measured by the number
of steps needed to get to the polytope by successively taking
pyramids over an initial polytope: the more steps we need the
simpler the polytope is. The following are immediately observed:

\smallskip$\centerdot$ the complexity is an invariant of the combinatorial type
and it never exceeds the dimension,

\smallskip$\centerdot$ a positive dimensional polytope $P$ is not a pyramid if and only if $\cc(P)=\dim P$,

\smallskip$\centerdot$ simplices are exactly the polytopes of complexity $0$,

\smallskip$\centerdot$ we always have the equality $\cc(\pyr(v,P))=\cc(P)$.

\

For a cone $C\subset\RR^n$ its \emph{complexity} $\cc(C)$ is
defined to be $\cc(\Phi(C))$ where $\Phi(C)=\G\cap C$ for any
affine hyperplane $\G\subset\RR^n$ cross-secting $C$. For a
positive affine monoid $M\subset\ZZ^n$ its \emph{complexity}
$\cc(M)$ is defined to be that of the cone $\RR_+M$.

Consider two polytopes $P\subset Q$, $P\not=Q$. Assume $P$ is
obtained from $Q$ by cutting off a pyramid at a vertex $v\in Q$.
In other words, $Q=P\cup P'$, $\dim P=\dim P'=\dim Q$ and
$P'=\pyr(v,P\cap P')$. In this situation we say that $P\subset Q$
is a \emph{pyramidal extension}. Observe that if $P\subset Q$ is a
pyramidal extension then $\dim P=\dim Q\ge1$.

\smallskip The following lemma is a key combinatorial fact. Let
$P\subset\RR^n$ be a polytope. Call a sequence of polytopes
$P=P_0,P_1,P_2,\ldots$ \emph{admissible} if the following
conditions hold for all indices $k$:

\smallskip$\centerdot$ either $P_{k+1}\subset P_k$ is a pyramidal
extension or $P_k\subset P_{k+1}$,

\smallskip$\centerdot$ $P_k\subset P$.

\smallskip(Observe, $\dim P_k=\dim P_0$ for all $k$.)

\begin{lemma}\label{admis}
Let $P$ be a polytope and $U\subset P$ an open subset. There
exists an admissible sequence of polytopes $P=P_0,P_1,P_2,\ldots$
such that $P_j\subset U$ for all sufficiently large $j$. If $P$ is
rational then the polytopes $P_j$ can be chosen to be rational.

\smallskip If $P$ and $U$ are given explicitly (say, by the
vertices or support hyperplanes of $P$ and of a simplex inside
$U$). Then there is an algorithm that finds an admissible sequence
$P=P_0,P_1,P_2,\ldots$
\end{lemma}

The lemma is proved in \cite{G1} without explicit reference to the
algorithmic aspect (see \cite[\S8.G]{BrG} for the most recent
exposition). However, the proof is in fact algorithmic, see for
instance \cite{LW}.

\subsection{Sufficiency of pyramidal descent}\label{pyrsuf-b}
An extension of monoids $L\subset N$ is called \emph{pyramidal}
if:

\smallskip$\centerdot$ $L,N\subset\ZZ^n$ are nonzero affine
positive normal monoids,

\smallskip$\centerdot$ $\Phi(L)\subset\Phi(N)$ is a pyramidal extension of
polytopes,

\smallskip$\centerdot$ $N|\Phi(L)=L$.

\noindent\smallskip Here $\Phi(N)=\RR_+N\cap\G$ and for an
arbitrarily fixed rational affine hyperplane $\G\subset\RR^n$
cross-secting the cone $\RR_+N$.

\smallskip Observe that if $L\subset N$ is a pyramidal extension
then $\rank L=\rank N\ge2$.

\smallskip Let $L\subset N$ be a pyramidal extension of monoids. It
will be called \emph{an extension of complexity} $c$ if
$\cc\bigl(\overline{\Phi(N)\setminus\Phi(L)}\bigr)=c$, where
$\overline{\mathcal Z}$ refers to the closure of $\mathcal Z$ in
the Euclidean topology. In this situation we will write
$\cc(L\subset N)=c$.

\smallskip We say that $\GL_r$-\emph{pyramidal descent} holds for a
pyramidal extension of monoids $L\subset N$ if for every
explicitly given matrix $A\in\GL_r(R[N_*])$ one can effectively
find a natural number $j$ and an elementary matrix
$E\in\E_r(R[(N_*)^{c^{-j}}])$, together with a representation
$E=\prod_ke_{p_kq_k}(\lambda_k)$ where
$e_{p_kq_k}(\lambda_k)\in\E_r(R[(N_*)^{c^{-j}}])$, such that
$EA\in\GL_r(R[(L_*)^{c^{-j}}])$. We say that
\emph{$\GL_r$-pyramidal descent of type $c$ holds for monoids of
rank $m$} for some $m\in\NN$ if $\GL_r$-pyramidal descent holds
for all pyramidal extensions of monoids $L\subset N$ with
$\cc(L\subset N)=c$ and $\rank N=m$.

\smallskip We say that $K_{2,r}$-pyramidal descent holds for a pyramidal
extension of monoids $L\subset N$ if the homomorphism
$K_{2,r}(R[(L_*)^{c^{-\infty}}])\to
K_{2,r}(R[(N_*)^{c^{-\infty}}])$ is surjective.  We say that
\emph{$K_{2,r}$-pyramidal descent of type $c$ holds for monoids of
rank $m$} for some $m\in\NN$ if pyramidal descent holds for all
pyramidal extensions $L\subset N$ with $\cc(L\subset N)=c$ and
$\rank N=m$.

\begin{proposition}\label{sufficient} Let $M\subset\ZZ^n$ be an affine positive
normal monoid. Then Theorem \ref{MAIN}(a) (corr. Theorem
\ref{MAIN}(b,c)) holds for the monoid algebra $R[M_*]$ if
$K_{2,r}$-pyramidal descent ($\GL_r$-pyramidal descent) of type
$<\cc(M)$ holds for monoids of rank $=\rank M$. (Here $r$ is as in
the corresponding part of Theorem \ref{MAIN}.)
\end{proposition}

\begin{proof}
Let $P_0\subset P_1\subset\cdots\subset P_i=\Phi(M)$ be a sequence
of rational polytopes where
$i=\rank(M)-1-\cc(M)=\dim\Phi(M)-\cc(M)$ and
$P_j=\pyr(v_j,P_{j-1})$ for each $j\in[1,i]$.

Fix a rational simplex $\Delta\subset P_0$, $\dim\Delta=\dim P_0$.
By Lemma \ref{admis} there is an admissible sequence
$P_0=Q_0,Q_1,Q_2,\ldots$ of rational polytopes such that
$Q_k\subset\Delta$ for all $k\gg0$. Then the sequence of polytopes
$\tilde Q_k=\conv(v_1,\ldots,v_i,Q_k)$ is an admissible sequence
of rational polytopes such that $\tilde Q_0=\Phi(M)$ and $\tilde
Q_k$ are contained in the simplex
$\tilde\Delta=\conv(v_1,\ldots,v_i,\Delta)$ for $k\gg0$. (We
assume $\tilde Q_k=Q_k$ and $\tilde\Delta=\Delta$ when $i=0$.)
Moreover, if $Q_{k+1}\subset Q_k$ is a pyramidal extension then we
have
$$
\cc\bigl(\overline{\tilde Q_{k+1}\setminus\tilde
Q_k}\bigr)\leq\cc(M)-1.
$$
By Gordan's lemma (see Section \ref{Monoids}) the monoids
$\RR_+\tilde Q_t\cap M$ are all affine. Obviously, they are also
normal and positive.

\smallskip By Corollary \ref{internormred}, for Theorem \ref{MAIN}(a) it
is enough to show that
$$
K_2(R)=K_{2,r}(R[(M_*)^{c^{-\infty}}]).
$$

Let $x\in K_{2,r}(R[(M_*)^{c^{-\infty}}])$. Assume
$K_{2,r}$-pyramidal descent of type $<\cc(M)$ holds for monoids of
rank $=\rank M$. Then there exist a sequence of elements
$$
x_k\in K_{2,r}(R[(M|\tilde Q_k)_*]),\quad k=0,1,\ldots
$$
such that:

$\centerdot$ $x_0=x$,

$\centerdot$ if $\tilde Q_{k+1}\subset\tilde Q_k$ is a pyramidal
extension for some $k\ge0$ then $x_k$ is the image of $x_{k+1}$
under the map
$$
K_{2,r}(R[((M|\tilde Q_{k+1})_*)^{c^{-\infty}}])\to
K_{2,r}(R[((M|\tilde Q_k)_*)^{c^{-\infty}}]),
$$

$\centerdot$ if $\tilde Q_k\subset\tilde Q_{k+1}$ then $x_{k+1}$
is the image of $x_k$ under the map
$$
K_{2,r}(R[((M|\tilde Q_k)_*)^{c^{-\infty}}])\to
K_{2,r}(R[((M|\tilde Q_{k+1})_*)^{c^{-\infty}}]).
$$

\noindent In particular, for $k\gg0$ we have
$$
x\in\Im(K_{2,r}(R[((M|\tilde\Delta)_*)^{c^{-\infty}}])\to
K_{2,r}(R[((M_*)^{c^{-\infty}}]).
$$
In view of Theorems \ref{BHS} and \ref{Tulenbaev} Proposition
\ref{apprA} we are done.

\

The case of Theorem \ref{MAIN}(b,c) is treated by the obvious
adaptation of the argument above, using the $\GL_r$-pyramidal
descent. For the algorithmic issues it is of course important that
all the involved convex polyhedral constructions can be carried
out effectively.
\end{proof}

In view of the equation (\ref{M*}) in Section \ref{seminormal} and
Proposition \ref{efectmiln} we get

\begin{corollary}\label{actuallysufficient}
Theorem \ref{MAIN} follows if $K_{2,r}$ and $\GL_r$-pyramidal
descents hold for any pyramidal extension of monoids, where $r$ is
as in the corresponding part of Theorem \ref{MAIN}.
\end{corollary}

\begin{remark}\label{whycompl}
As mentioned, the concept of a pyramidal descent without
consideration of complexities was introduced in \cite{G1}: using
induction on $\rank N$, it is shown in \cite{G1} that (unstable)
$K_0$-pyramidal descent holds for all pyramidal extensions
$L\subset N$. The complexities were added to the picture in
\cite{G5} for reasons not related to this paper at all. However,
it is the notion of complexity that makes the induction argument
work in Section \ref{PYRAMID} where we show that, indeed, $\GL_r$
and $K_{2,r}$-pyramidal descents hold for all pyramidal extensions
$L\subset N$. The argument will use induction on the pairs $(\rank
N,\cc(L\subset N))$. In \cite{Mu} this aspect is simply absent.
\end{remark}

\section{Almost separation}\label{ALMOST}
In this section we state the main technical fact to be used in the
proof of $K_{2,r}$- and $\GL_r$-pyramidal decents.

\smallskip Let $M\subset\ZZ^n$ be an affine positive normal monoid with
$\gp(M)=\ZZ^n$.

\smallskip Let $\H\subset\RR^n$ be a rational hyperplane, dissecting
the cone $\RR_+M$ into two $n$-cones $\RR_+M=C_1\cup C_2$. Fix a
rational affine hyperplane $\G\subset\RR^n$ with
$\RR_+M=\RR_+(\RR_+M\cap\G)$.

We also fix a real number $\epsilon>0$ and a natural number
$c\ge2$.

Let $M_1(\epsilon)=\RR_+M\cap C_1(\epsilon)\cap M$ and
$M_2(\epsilon)=\RR_+M\cap C_2(\epsilon)\cap M$, where
$C_1(\epsilon)$ and $C_2(\epsilon)$ refer to the open cones
introduced in Section \ref{polytopesandcones}.

For a ring $\Lambda$ and a matrix $A\in\E_r(\Lambda)$ under
\emph{a representation $\bar A$} we will mean a representation of
the form
$$
A=\prod_ke_{p_kq_k}(\lambda_k),\quad \lambda_k\in\Lambda.
$$

The theorem below is essentially due to Mushkudiani \cite{Mu} (see
Remark \ref{mushkwrong}):

\begin{theorem}\label{elparation}
Let $R$ be an arbitrary ring, $r\ge2$ be a natural number and
$A\in\E_r(R[M_*])$. Then for any representation $\bar A$ one can
explicitly find a natural number $j_A$ and a factorization of the
form $A=A_1 A_2$ for some
$A_1\in\E_r(R[(M_1(\epsilon)_*)^{c^{-j_A}}])$, together with a
representation $\bar A_1$, and
$A_2\in\SL_r(R[(M_2(\epsilon)_*)^{c^{-j_A}}])$.
\end{theorem}

\smallskip\noindent(The equality $A=A_1A_2$ is considered
in the ambient group $\GL_r(R[M^{c^{-\infty}}])$.)

\smallskip In other words, the input of the algorithm is an
explicit representation of the form $\bar
A=\prod_ke_{p_kq_k}(\lambda_k)$, $\lambda_k\in R[M_*]$, and the
output is a natural number $j_A$ and a factorization $A=A_1A_2$
where $A_1\in\E_r(R[(M_1(\epsilon)_*)^{c^{-j_A}}])$ and
$A_2\in\SL_r(R[(M_2(\epsilon)_*)^{c^{-j_A}}])$, together with an
explicit representation of the form
$$
\bar A_1=\prod_ke_{r_ks_k}(\mu_k),\quad \mu_k\in
R[(M_1(\epsilon)_*)^{c^{-j_A}}].
$$
Here it is assumed that in $R$ and $M$ we can explicitly perform
the operations.

We want to emphasize that even \emph{without referring the
algorithmic aspect}, Theorem \ref{elparation} states a nontrivial
fact which leads to the nilpotence of $K_{1,r}(R[M])$.

\begin{remark}\label{thename}
In view of Theorem \ref{MAIN}(b), Theorem \ref{elparation} is
equivalent to the equality
$$
\E_r(R[(M_*)^{c^{-\infty}}])=\E_r(R[(M_1(\epsilon)_*)^{c^{-\infty}}])
\E_r(R[(M_2(\epsilon)_*)^{c^{-\infty}}])
$$
which actually explains the name `almost separation'. However,
since Theorem \ref{MAIN} is a consequence of Theorem
\ref{elparation}, we have to resort to the formulation above.
\end{remark}

\begin{remark}\label{mushkwrong}
Mushkudiani's original version, derived in the course of the proof
of \cite[Theorem 3.1]{Mu} (but not stated explicitly), claims the
existence of a representation of the form $A=A_1A_2$ where
$A_1\in\E(R[(\RR_+M\cap C_1\cap M)^{c^{-\infty}}])$ and
$A_2\in\SL(R[M_2(\epsilon)^{c^{-\infty}}])$. However, the
corrected argument, presented in Section \ref{Mushk}, gives the
current version. Moreover, the argument in \cite{Mu} never really
uses the fact that in Theorem \ref{elparation} one takes iterated
$c$th roots of monomials. But without taking the $c$th roots of
monomials, Theorem \ref{elparation} can \emph{not} hold as it
would lead to a contradiction with \cite{G3} and \cite{Sr}.
\end{remark}

The next theorem is a $\st$-version of Theorem \ref{elparation}.

\begin{theorem}\label{steparation}
Let $R$ be a regular ring and $r\ge\max(5,\dim R+3)$ be a natural
number. Then any element $x\in\st_r(R[(M_*)^{c^{-\infty}}])$ has a
factorization of the form:
\begin{align*}
x=yz,\qquad &y\in\Im(\st_r(R[(M_1(\epsilon)_*)^{c^{-\infty}}])\to
\st_r(R[(M_*)^{c^{-\infty}}])),\\
&z\in\Im(\st_r(R[(M_2(\epsilon)_*)^{c^{-\infty}}])
\to\st_r(R[(M_*)^{c^{-\infty}}])).
\end{align*}
\end{theorem}

The logical scheme of the relationships between Theorems
\ref{MAIN}, \ref{elparation} and \ref{steparation} is given by the
following diagram:
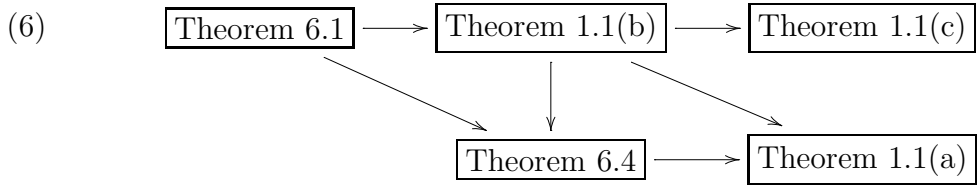
\begin{equation}\label{logics}
\xymatrix{\fbox{\text{Theorem
\ref{elparation}}}\ar[r]\ar[rd]&\fbox{\text{Theorem
\ref{MAIN}(b)}}\ar[r]\ar[d]\ar[rd]&\fbox{{\text{Theorem \ref{MAIN}(c)}}}\\
&\fbox{\text{Theorem
\ref{steparation}}}\ar[r]&\fbox{{\text{Theorem \ref{MAIN}(a)}}}}
\end{equation}
which will be realized gradually in the following sections,
postponing the proof of Theorem \ref{elparation} to the very end.

\smallskip Below we explain how Theorems \ref{MAIN}(b) and Theorem \ref{elparation}
together imply Theorem \ref{steparation}. This corresponds to the
left triangle in diagram (\ref{logics}).

\begin{proof} For simplicity of notation let
\begin{align*}
&\mathcal Y=\Im(\st_r(R[(M_1(\epsilon)_*)^{c^{-\infty}}])\to
\st_r(R[(M_*)^{c^{-\infty}}])),\\
&\mathcal Z=\Im(\st_r(R[(M_2(\epsilon)_*)^{c^{-\infty}}])
\to\st_r(R[(M_*)^{c^{-\infty}}])).
\end{align*}

\smallskip First we consider the case when $M$ is simplicial.

Let $E\in\E_r(R[(M_*)^{c^{-\infty}}])$ denote the image of $x$. By
Theorem \ref{elparation} we can write $E=A_1A_2$ where
$A_1\in\E_r(R[(M_1(\epsilon)_*)^{c^{-\infty}}])$ and
$A_2\in\SL_r(R[(M_2(\epsilon)_*)^{c^{-\infty}}])$. By Theorem
\ref{MAIN}(b) (or, equivalently, Theorem \ref{cdivrev}(b))
$A_2\in\E_r(R[(M_2(\epsilon)_*)^{c^{-\infty}}])\SL_r(R)$.

Since $R$ is a retract of the rings $R[(M_*)^{c^{-\infty}}]$,
$R[(M_1(\epsilon)_*)^{c^{-\infty}}]$ and
$R[(M_2(\epsilon)_*)^{c^{-\infty}}]$, we actually have
$A_2\in\E_r(R[(M_2(\epsilon)_*)^{c^{-\infty}}])$.

By lifting $A_1$ and $A_2$ respectively to $\mathcal Y$ and
$\mathcal Z$ we find two elements $y\in\mathcal Y$ and
$z\in\mathcal Z$ such that $x=\xi yz$ for some $\xi\in
K_{2,r}(R[(M_*)^{c^{-\infty}}])$. By Proposition \ref{apprA} the
monoid $(M_*)^{c^{-\infty}}$ is a filtered union of free monoids.
Therefore, Theorems \ref{BHS} and \ref{Tulenbaev} imply that
$K_{2,r}(R[(M_*)^{c^{-\infty}}])=K_2(R)$. In particular, $\xi
z\in\mathcal Y$. Hence the desired representation $x=(\xi y)z$.

\

Now we consider the case of a general affine positive normal
monoid $M\subset\ZZ^n$.

\smallskip Fix a surjective monoid homomorphism $\pi:\ZZ_+^m\to M$
for some $m$. Its $\RR$-linear extension $\RR^m\to\RR^n$ will be
denoted by $\RR\otimes\pi$.

There exist a rational hyperplane $\mathcal H'\subset\RR^m$,
dissecting the standard positive orthant $\RR_+^m$ into two
$m$-cones $\RR_+^m=C_1'\cup C_2'$, and a real number $\epsilon'>0$
such that $(\RR\otimes\pi)(C_1'(\epsilon'))\subset C_1(\epsilon)$
and $(\RR\otimes\pi)(C_2'(\epsilon'))\subset C_2(\epsilon)$. Here
the open convex cones $C_1'(\epsilon')$ and $C_2'(\epsilon')$ are
considered with respect to arbitrarily fixed affine hyperplane
$\G'\subset\RR^m$, cross-secting the positive orthant $\RR_+^m$.
Let $\mathcal Y'$ and $\mathcal Z'$ denote the sets
\begin{align*}
&\mathcal Y'=\Im(\st_r(R[(M_1'(\epsilon')_*)^{c^{-\infty}}])\to
\st_r(R[((\ZZ_+^m)_*)^{c^{-\infty}}])),\\
&\mathcal Z'=\Im(\st_r(R[(M_2'(\epsilon')_*)^{c^{-\infty}}])
\to\st_r(R[((\ZZ_+^m)_*)^{c^{-\infty}}])).
\end{align*}
Then $\pi$ induces a surjective group homomorphism
$$
\pi_*:\st_r(R[((\ZZ_+^m)_*)^{c^{-\infty}}]\to
\st_r(R[(M_*)^{c^{-\infty}}])
$$
such that $\pi_*(\mathcal Y')\subset\mathcal Y$ and
$\pi_*(\mathcal Z')\subset\mathcal Z$. Therefore, the general case
reduces to the case when $M$ is simplicial.
\end{proof}

\begin{remark}\label{mushksteparation}
The proof of Theorem \ref{steparation} (in a slightly different
formulation) constitutes the main part of \cite{Mu}. It represents
a `Steinberg group version' of the argument in Section
\ref{Mushk}. However, the approach in \cite{Mu} simply cannot be
rescued.
\end{remark}

\begin{remark}\label{gencuts}
As it becomes clear in Section \ref{proofmain}, we only need the
validity of Theorems \ref{elparation} and \ref{steparation} for
the special cuts of $\RR_+M$ by $\mathcal H$ when one extremal ray
of $\RR_+M$ lies strictly on one side of $\mathcal H$ and the
other extremal rays lie on the other side. However, our deduction
of Theorem \ref{steparation} from Theorem \ref{elparation} is
through lifting the general case to the case when $M$ is
simplicial (the map $\pi$ above) and the mentioned condition on
the dissecting hyperplane is in general \emph{not} respected under
such a lifting. So we really need the general version of Theorem
\ref{elparation}.
\end{remark}

\section{Almost separation implies pyramidal descent}\label{PYRAMID}
In this section $R$ is a regular ring of finite Krull dimension.

\smallskip In Section \ref{mainb} we assume $r\ge\max(3,\dim
R+2)$ and show how Theorem \ref{elparation} implies
\emph{$\GL_r$-pyramidal descent}. This corresponds to the upper
left horizontal arrow in diagram (\ref{logics}). The upper right
arrow simply reflects the fact that the proof of Theorem
\ref{MAIN}(b) is algorithmic in nature.

\smallskip In Section \ref{maina} we assume $r\ge\max(5,\dim
R+3)$ and show how Theorems \ref{MAIN}(b) and \ref{steparation}
imply \emph{$K_{2,r}$-pyramidal descent}. This corresponds to the
right triangle in diagram (\ref{logics}).

\subsection{$\GL_r$-pyramidal descent}\label{mainb} Here we prove

\begin{lemma}\label{GLrpyrdescent}
$GL_r$-pyramidal descent holds for any pyramidal extension of
monoids.
\end{lemma}

\begin{proof} Let  $L\subset N$ be a pyramidal extension of monoids in $\ZZ^n$.
We use induction on the pairs $(\rank N, \cc(N))$, ordered
lexicographically.

If $\cc(N)=0$ then $N$ is simplicial and then we are done by
Corollaries \ref{simplicialalg} and \ref{simplicialK1}. Notice,
the condition $\cc(N)=0$ also includes the case $\rank N\le2$.

Now assume $\cc(N)>0$ and the $\GL_r$-pyramidal descent has been
shown for the pyramidal extensions $L'\subset N'$ for which
$$
(\rank N',\cc(L'\subset N'))<(\rank N,\cc(L\subset N)).
$$
We want to show the equality
\begin{equation}\label{WTS}
\GL_r(R[(N_*)^{c^{-\infty}}])=\E_r(R[(N_*)^{c^{-\infty}}])\GL_r(R[(L_*)^{c^{-\infty}}]).
\end{equation}

\smallskip By Proposition \ref{sufficient} for any affine positive normal monoid
$M'$, satisfying the conditions $ \rank M'=\rank N$ and
$\cc(M')=\cc(L\subset N)$, we have
\begin{equation}\label{M'*}
\GL_r(R[(M'_*)^{c^{-\infty}}])=\E_r(R[(M'_*)^{c^{-\infty}}])\GL_r(R).
\end{equation}

\smallskip Fix an affine hyperplane $\G\subset\RR^n$ cross-secting
the cone $\RR_+N$. The $\Phi$-polytopes below are all considered
w.r.t. $\G$.

$\Phi(N)$ has exactly one vertex that does not belong to
$\Phi(L)$. Call it $v$. Let $C(v,\Phi(N))\subset\G$ denote the
affine cone spanned by $\Phi(N)$ at $v$, that is
$$
C(v,\Phi(N))=v+\RR_+\left(\Phi(N)-v\right).
$$

We have the rational pyramid
$\Delta_1=\overline{\Phi(N)\setminus\Phi(L)}\subset\Phi(N)$.

Let $\Delta_2\subset C(v,\Phi(N))$ be any rational pyramid
satisfying the conditions:

\smallskip$\centerdot$ $v\in\vert(\Delta_2)$,

\smallskip$\centerdot$ $C(v,\Phi(N))=v+\RR_+\left(\Delta_2-v\right)$,

\smallskip$\centerdot$ $\Phi(N)\subset\Delta_2$.

\

\noindent The following two conditions are satisfied
automatically:

\smallskip$\centerdot$  $\cc(\Delta_2)=\cc(\Delta_1)=\cc(L\subset
N)$,

\smallskip$\centerdot$ $\dim\Delta_2=\dim\Delta_1=\dim\Phi(N)$.

\

\noindent In particular, (\ref{M'*}) implies
\begin{equation}\label{indhypothesis}
\GL_r(R[(N_*)^{c^{-\infty}}])\subset\GL_r(R[((N|\Delta_2)_*)^{c^{-\infty}}])
=\E_r(R[((N|\Delta_2)_*)^{c^{-\infty}}])\GL_r(R).\footnote{It is
here where we need $N$ to be normal -- it enables us to consider
the monoid $N|\Delta_2$ which satisfies the condition
 $(N|\Delta_2)|\Phi(N)=N$.}
\end{equation}

\smallskip Fix a rational point $\xi\in\inte(\Phi(L))$.
For a real number $\lambda$ the homothetic image of a polytope
$\Pi\subset\G$ with the factor $\lambda\in\RR$ and centered at
$\xi$ will be denoted by $\Pi_\lambda$.

\smallskip For any real number $0<\lambda<1$ we fix a real number
$\epsilon_\lambda>0$ in such a way that

\begin{equation}\label{second}
\Phi(L)_\lambda(\epsilon_\lambda)\subset\inte(\Phi(L)).
\end{equation}

\smallskip  Furthermore, for a rational number $0<\lambda<1$ we use the notation:
\begin{align*}
N_{1,\lambda}(\epsilon_\lambda)_*=
((N|(\Delta_1)_\lambda(\epsilon_\lambda))_*)^{c^{-\infty}}\quad\text{and}\quad
N_{2,\lambda}(\epsilon_\lambda)_*=((N|(\overline{\Delta_2\setminus\Delta_1})_\lambda
(\epsilon_\lambda))_*)^{c^{-\infty}},
\end{align*}
where $(\Delta_1)_\lambda(\epsilon_\lambda)$ and
$(\overline{\Delta_2\setminus\Delta_1})_\lambda(\epsilon_\lambda)$
correspondingly refer to the $\epsilon_\lambda$-neighborhoods of
$(\Delta_1)_\lambda$ and
$(\overline{\Delta_2\setminus\Delta_1})_\lambda$ inside the
pyramid $(\Delta_2)_\lambda$.

We record the following consequence of (\ref{second}):
\begin{equation}\label{third}
\inte(\Phi(N))\cap(\overline{\Delta_2\setminus\Delta_1})_\lambda(\epsilon_\lambda)\subset\inte(\Phi(L)).
\end{equation}
((\ref{second}) guarantees that the part of
$(\overline{\Delta_2\setminus\Delta_1})_\lambda(\epsilon_\lambda)$
`towards $v$' is in $\inte(\Phi(L))$.)

\smallskip Now by Theorem \ref{elparation}  we have
$$
\E_r(R[((N|(\Delta_2)_\lambda)_*)^{c^{-\infty}}])\subset
\E_r(R[N_{1,\lambda}(\epsilon_\lambda)_*])\SL_r(R[N_{2,\lambda}(\epsilon_\lambda)_*])
$$
which, in view of (\ref{indhypothesis}), implies
\begin{equation}\label{tobeunited}
\GL_r(R[((N|(\Delta_2)_\lambda)_*)^{c^{-\infty}}])\subset
\E_r(R[(N_*)^{c{^{-\infty}}}])\GL_r(R[N_{2,\lambda}(\epsilon_\lambda)_*]).
\end{equation}
By letting $\lambda$ run over the set $\QQ\cap(0,1)$, the
inclusion (\ref{tobeunited}) implies
\begin{equation}\label{thirteen}
\begin{aligned}
\GL_r(R[(N_*)^{c^{-\infty}}])\subset \bigcup_\lambda\
\E_r(R[(N_*)^{c{^{-\infty}}}])\GL_r(R[N_{2,\lambda}(\epsilon_\lambda)_*]).
\end{aligned}
\end{equation}

Now (\ref{WTS}) follows from (\ref{thirteen}) once we show the
following implication for any $\lambda$:
\begin{align*}
\begin{cases}
A=BC\\
A\in\GL_r(R[(N_*)^{c^{-\infty}}])\\
B\in\E_r(R[(N_*)^{c{^{-\infty}}}])\\
C\in\GL_r(R[N_{2,\lambda}(\epsilon_\lambda)_*])
\end{cases}\quad
\Longrightarrow\quad C\in\GL_r(R[(L_*)^{c^{-\infty}}]).
\end{align*}
But for such a triple of matrices, using (\ref{third}), we have
\begin{align*}
C=&B^{-1}A\in\GL_r(R[(N_*)^{c^{-\infty}}])\cap\GL_r(R[N_{2,\lambda}(\epsilon)_*])=\\
&\GL_r(R[(N_*)^{c^{-\infty}}\cap
N_{2,\lambda}(\epsilon)_*])=\GL_r(R[(N|\inte(\Phi(N))\cap
(\overline{\Delta_2\setminus\Delta_1})_\lambda(\epsilon_\lambda))^{c^{-\infty}}])\subset\\
&\GL_r(R[(N|\inte(\Phi(L)))^{c^{-\infty}}])=\GL_r(R[(L_*)^{c^{-\infty}}]).
\end{align*}
\end{proof}

\subsection{$K_{2,r}$-pyramidal descent}\label{maina} Here we
prove
\begin{lemma}\label{k2rpydescent}
$K_{2,r}$-pyramidal descent holds for any pyramidal extension of
monoids $L\subset N$.
\end{lemma}

\begin{proof}
We use the same induction as in the proof of Lemma
\ref{GLrpyrdescent}, that is the induction on the pairs $(\rank N,
\cc(N))$, ordered lexicographically. Also, we assume that
$L,N\subset\ZZ^n$.

If $\cc(N)=0$ then $N$ is simplicial and then we are done by
Corollary \ref{simplK2}. This also includes the case $\rank
N\le2$.

Now assume $\cc(N)>0$ and $K_{2,r}$-pyramidal descent has been
shown for the pyramidal extensions $L'\subset N'$ for which
$$
(\rank N',\cc(L'\subset N'))<(\rank N,\cc(L\subset N)).
$$
Pick an arbitrary element $x\in K_{2,r}(R[(N_*)^{c^{-\infty}}])$.
We want to show
\begin{equation}\label{WTS2}
x\in\Im\big(K_{2,r}(R[(L_*)^{c^{-\infty}}])\to
K_{2,r}(R[(N_*)^{c^{-\infty}}])\big).
\end{equation}

\smallskip By Proposition \ref{sufficient} for any affine positive normal monoid
$M'$, satisfying the conditions $ \rank M'=\rank N$ and
$\cc(M')=\cc(L\subset N)$, we have
\begin{equation}\label{M'*2}
K_2(R)=K_{2,r}(R[(M'_*)^{c^{-\infty}}]).
\end{equation}

Fix a rational affine hyperplane $\G\subset\RR^n$ cross-secting
the cone $\RR_+N$. The $\Phi$-polytopes below are all considered
w.r.t. $\G$.

We have the pyramid $\Delta=\overline{\Phi(N)\setminus\Phi(L)}$.
Fix a rational point $\xi\in\inte(\Phi(L))$, a rational number
$0<\lambda<1$ and a real number $\epsilon>0$ so that the following
conditions are satisfied\footnote{This can be done first by
choosing $\lambda$ sufficiently close to 1 and then choosing
$\epsilon$ sufficiently small, depending on $\lambda$.}:

\smallskip$\centerdot$ $x$ is the image of some $x_\lambda\in
K_{2,r}(R[((N_\lambda)_*)^{-\infty}])$ where
$N_\lambda=N|\Phi(N)_\lambda$,

\smallskip$\centerdot$
$\Phi(L)_\lambda(\epsilon)\subset\inte(\Phi(L))$

\smallskip$\centerdot$ $\Delta_\lambda(\epsilon)\subset\Delta'$
for some rational simplex $\Delta'\subset\inte(\Phi(N))$, similar
to $\Delta$.

\smallskip\noindent Above we have used the notation:

\smallskip$\centerdot$ for any polytope $\Pi\subset\Phi(N)$ its homothetic image with
factor $\lambda$ and centered at $\xi$ is denoted by
$\Pi_\lambda$,

\smallskip$\centerdot$ for any polytope $\Pi\subset\Phi(N)$
its $\epsilon$-neighborhood inside $\Phi(N)$ is denoted by
$\Pi(\epsilon)$.

\

Consider the monoids
$M_1(\epsilon)=N_\lambda|\Delta_\lambda(\epsilon)$ and
$M_2(\epsilon)=N_\lambda|\Phi(L)_\lambda(\epsilon)\subset L_*$. By
Theorem \ref{steparation} we have a representation of the form:
\begin{align*}
&x_\lambda=yz,\\
&y\in\Im(\st_r(R[(M_1(\epsilon)_*)^{c^{-\infty}}])\to
\st_r(R[((N_\lambda)_*)^{c^{-\infty}}])),\\
&z\in\Im(\st_r(R[(M_2(\epsilon)_*)^{c^{-\infty}}])\to
\st_r(R[((N_\lambda)_*)^{c^{-\infty}}])).
\end{align*}
For the corresponding elementary matrices
$E_y,E_z\in\E_r(R[(N_*)^{c^{-\infty}}])$ we have
\begin{align*}
E_yE_z=1,\quad
E_y\in\E_r(R[(M_1(\epsilon)_*)^{c^{-\infty}}]),\quad
E_z\in\E_r(R[(M_2(\epsilon)_*)^{c^{-\infty}}]),
\end{align*}
which implies
$$
E_y,E_z\in\SL_r(R[(M_1(\epsilon)_*)^{c^{-\infty}}\cap
(M_2(\epsilon)_*)^{c^{-\infty}}])=\SL_r(R[(M_1(\epsilon)_*\cap
M_2(\epsilon)_*)^{c^{-\infty}}])
$$
By Theorem \ref{MAIN}(b) we get
$$
E_y,E_z\in\E_r(R[(M_1(\epsilon)_*\cap
M_2(\epsilon)_*)^{c^{-\infty}}]).
$$
Let
$$
w\in\Im\big(\st_r(R[(M_1(\epsilon)_*\cap
M_2(\epsilon))_*^{c^{-\infty}}])\to\st_r(R[((N_\lambda)_*)^{c^{-\infty}}])\big)
$$
be any lifting of $E_y$. Then we have:
\begin{align*}
&x_\lambda=(yw^{-1})\cdot(wz),\\
&yw^{-1}\in\Im(\st_r(R[(M_1(\epsilon)_*)^{c^{-\infty}}])\to
\st_r(R[((N_\lambda)_*)^{c^{-\infty}}])),\\
&wz\in\Im(\st_r(R[(M_2(\epsilon)_*)^{c^{-\infty}}])\to
\st_r(R[((N_\lambda)_*)^{c^{-\infty}}])),
\end{align*}

Since the image of $yw^{-1}$ in
$\E_r(R[(M_1(\epsilon)_*)^{c^{-\infty}}])$ is $\1$ we actually
have
$$
yw^{-1}\in\Im(K_{2,r}(R[(M_1(\epsilon)_*)^{c^{-\infty}}])\to
K_{2,r}(R[((N_\lambda)_*)^{c^{-\infty}}]))
$$
and, similarly,
$$
wz\in\Im(K_{2,r}(R[(M_2(\epsilon)_*)^{c^{-\infty}}])\to
K_{2,r}(R[((N_\lambda)_*)^{c^{-\infty}}])).
$$
But then the inclusion $M_2(\epsilon)_*\subset L_*$ implies
$$
wz\in\Im(K_{2,r}(R[(L_*)^{c^{-\infty}}])\to
K_{2,r}(R[(N_\lambda)^{c^{-\infty}}])).
$$
In particular, (\ref{WTS2}) follows if we show that the image of
$yw^{-1}$ in $K_{2,r}(R[(N_*)^{c^{-\infty}}])$ belongs to
$K_2(R)$.

We have
\begin{align*}
&\Im\big(K_{2,r}(R[(M_1(\epsilon)_*)^{c^{-\infty}}])\to
K_{2,r}(R[(N_*)^{c^{-\infty}}])\big)\subset\\
&\Im\big(K_{2,r}(R[((N|\Delta')_*)^{c^{-\infty}}])\to
K_{2,r}(R[(N_*)^{c^{-\infty}}])\big)
\end{align*}
and, in view of the conditions $\rank(N|\Delta')=\rank N$ and
$\cc(N|\Delta')=\cc(L\subset N)$, by (\ref{M'*2}) we get
$K_{2,r}(R[((N|\Delta')_*)^{c^{-\infty}}])=K_2(R)$.
\end{proof}

\section{Proof of Theorem \ref{elparation}}\label{Mushk}

This section presents a corrected version of Mushkudiani's proof
of almost separation in $\E_r(R[M])$. The algorithmic part of
Theorem \ref{elparation} is a direct consequence of the argument
presented below and we do not discuss it separately.

\subsection{Convention and notation}\label{Notation}
Here we introduce the notation to be used in the rest of Section
\ref{Mushk}.

\

\noindent\underline{\emph{Monoids and cones}}. We fix an affine
positive monoid $M\subset\QQ^n$, $n=\rank\gp(M)\ge2$. We don't
require that $M$ is normal or $M\subset\ZZ^n$. Let
$M^+=M\setminus\{0\}$.

\smallskip For a point $z\in\RR^n$ its $n$th coordinate will be
denoted by $z_n$.

\smallskip Assume a rational hyperplane $\H\subset\RR^n$ cuts
$\RR_+M$ into two $n$-dimensional subcones. Without loss of
generality we will assume $\H=\RR^{n-1}\oplus0\subset\RR^n$ -- a
condition that can be achieved by a rational coordinate change.

\smallskip We can additionally assume that the cone $\RR_+M$ is
`acute' enough to have the following condition satisfied:
\begin{equation}\label{Acute}
\forall u,v\in\RR_+M\setminus\{0\}\qquad  \|u\|,\|v\|<\|u+v\|.
\end{equation}
In fact, without loss of generality we can assume that no negative
multiple of $e_1$ belongs to $M$ and then (\ref{Acute}) can be
achieved by applying to $M$ a linear transformation of the form
$e_1\mapsto e_1$ and $e_i\mapsto e_i+ke_1$ with $k\gg0$ for
$i\not=1$. Here $\{e_1,\ldots,e_n\}$ is the standard basis of
$\RR^n$.

\smallskip We also fix a rational affine
hyperplane $\G\subset\RR^n$ such that
$\RR_+M=\RR_+\big(\RR_+M\cap\G\big)$. Thus $\Phi(M)=\RR_+M\cap\G$.
Recall, for any submonoid $N\subset M$ we put
$\Phi(N)=\RR_+N\cap\G$.

\

\noindent\underline{\emph{Monomials}}. Let $R$ be a ring.
\emph{Monomials} in $R[M]$ are simply the elements of $M$.

\smallskip The products $a\mu\in R[M]$, $a\in R$, $\mu\in M$ are
\emph{terms}. If $a\not=0$ then $\mu$ is called \emph{the support
monomial of $a\mu$}. For a nonzero element $\gamma\in R[M]$ the
support monomials in the canonical expansion of $\gamma$ as a sum
of terms constitute the set of \emph{the support monomials of
$\gamma$}. It is denoted by $\supp(\gamma)$.

\smallskip For a nonzero term $z=a\mu\in R[M]$, $a\in R$, $\mu\in
M$, its \emph{length} $\|a\mu\|$ is just the Euclidean norm
$\|\mu\|$ in $\RR^n$. Let $z_n=\mu_n$.

\smallskip For a subset $I\subset\RR$ we put
$$
R[M]_I=\{\gamma\in R[M]\ |\ \mu_n\in I\quad\text{for every}\quad
\mu\in\supp(\gamma)\}\subset R[M].
$$
Thus $0\in R[M]_I$ for any subset $I\subset\RR$ and $R\subset
R[M]_I$ if $0\in I$.

\smallskip For a nonzero term $z=a\mu\in R[M]$, $a\in R$, $\mu\in
M$, and a nonzero element $\gamma\in R[M]$ we put
$\Phi(z)=\G\cap\RR_+\mu$ and $\Phi(\gamma)=\conv\{\Phi(\mu)\ |\
\mu\in\supp(\gamma)\}$. By convention, $\Phi(0)=\emptyset$. In
particular, $\Phi(\gamma)$ is always a polytope inside $\Phi(M)$.

\smallskip For an element $\gamma\in R[M]$ we say that $\gamma_n$ (or
$\Phi(\gamma)_n$, or $\|\gamma\|$) satisfies certain inequality if
the $n$th coordinate (respectively, the $n$th coordinate of the
$\Phi$-image, the length) of every element $\mu\in\supp(\gamma)$
satisfies the same inequality.

\smallskip For real numbers $l>0$ and $\epsilon$ consider the
subset
$$
\B'(\epsilon,l)=\{\gamma\in R[M]\ |\ l\le\|\gamma\|,\
\epsilon\le\Phi(\gamma)_n\}\subset RM^+.
$$

\

\noindent\underline{\emph{Matrices}}. Fix a natural number
$r\ge2$. For a matrix $A\in\M_r(R[M])$ a \emph{support monomial of
$A$} is by definition a support monomial of some entry of $A$. The
set of support monomials of $A$ is denoted by $\supp(A)$.

\smallskip For a matrix $A=(\lambda_{ij})_{i,j=1}^r\in\M_r(R[M])$ we say that
$A_n$ satisfies certain inequality if every $(\lambda_{ij})_n$
does so.

\smallskip For real numbers $l>0$ and $\epsilon$ we introduce the
following subsets of $\M_r(R[M])$:
\begin{align*}
&\A(\epsilon)=\{A\in\M_r(R[M])\ |\ 1\notin\supp(A)\ \text{and}\
\epsilon\le A_n\},\\
&\B(\epsilon,l)=\B'(\epsilon,l)^{r\times r},\\
&\D=\{D\in\M_r(R[M])\ |\ 1\notin\supp(D),\ D\ \text{is diagonal}\
\text{and}\
D_n\ge0\},\\
&\D_{>0}=\{D\in\M_r(R[M])\ | D\ \text{is diagonal}\ \text{and}\
D_n>0\}.
\end{align*}
Observe that all these matrices have entries from $RM^+$ and that
the zero matrix belongs to each of the mentioned classes of
matrices.

\smallskip As in the previous sections, \emph{a representation
$\bar E$} for a matrix $E\in \E_r(R[M])$ means a representation of
the form $E=\prod e_{ij}(\gamma_{ij})$, $\gamma_{ij}\in R[M]$.
Moreover, we say that $\bar E_n$ (resp. $\Phi(\bar E)_n$)
satisfies certain inequality if every $(\gamma_{ij})_n$ (resp.
$\Phi(\gamma_{ij})_n$) does so.

\subsection{Commuting rules for elementary
matrices}\label{Commuting}
\begin{lemma}\label{Com1}
Let $\epsilon_1$, $\epsilon$, $l$ be positive real numbers,
$i\not=j$ natural numbers, $D\in\D$, and $\alpha,\beta\in R[M]$
nonzero terms. Assume $|\alpha_n|<\epsilon_1\le\beta_n$. Then:
$$
\left(e_{ji}(\beta)+D\right)e_{ij}(\alpha)=e_{ij}(\alpha)e_{ij}(\gamma)
\left(\1+A+B+D'\right)
$$
for some
$$
\gamma\in
R[M]_{[\alpha_n,\epsilon_1)},\quad A\in\A(\epsilon_1),\quad
B\in\B(-\epsilon,l),\quad D'\in\D.
$$
Moreover, the support monomials of $\gamma$, $A$, $B$ and $D'$ are
products of those of $\alpha$, $\beta$ and $D$.
\end{lemma}

\smallskip\noindent(In this lemma we don't exclude the case $\alpha\in R$.)

\begin{proof} We want to find
$\gamma\in R[M]_{[\alpha_n,\epsilon_1)}$ and matrices $A,B,D'$ as
in the lemma such that
$$
e_{ij}(-\gamma)e_{ij}(-\alpha)\left(e_{ji}(\beta)+D\right)e_{ij}(\alpha)=
\1+A+B+D'.
$$
We have representations of the form:
\begin{itemize}
\item[$\centerdot$]
$e_{ij}(-\alpha)e_{ji}(\beta)e_{ij}(\alpha)=e_{ij}(a_0)+a_{ji}(\beta)+D_1$
for some  $a_0=-\alpha^2\beta\in R[M]_{(\alpha_n,+\infty)}$ and
$\D_1\in\D$, \item[$\centerdot$]
$e_{ij}(-\alpha)De_{ij}(\alpha)=D+a_{ij}(b_0)$ for some $b_0\in
R[M]_{[\alpha_n,+\infty)}$, \item[$\centerdot$]
$a_0+b_0=\gamma_1+a_1+b_1$ for some $\gamma_1\in
R[M]_{[\alpha_n,\epsilon_1)}\setminus \B'(-\epsilon,l)$, $a_1\in
R[M]_{[\epsilon_1,+\infty)}$, and $b_1\in\B'(-\epsilon,l)$.
\end{itemize}

\smallskip\noindent(Such a representation $a_0+b_0=\gamma_1+a_1+b_1$ is in general not
unique.)

\smallskip If $\gamma_1=0$ then we are done because
\begin{align*}
e_{ij}(-\alpha)\left(e_{ji}(\beta)+D\right)e_{ij}(\alpha)=
\1+\big(a_{ij}(a_1)+a_{ji}(\beta)\big)+a_{ij}(b_1)+(D+D_1).
\end{align*}

\smallskip So we can assume $\gamma_1\not=0$. Then we have
representations of the form:
\begin{align*}
&e_{ij}(-\gamma_1)\left(e_{ij}(\gamma_1)+D+D_1\right)=
e_{ij}(\delta_1)+D+D_1,
\quad\delta_1\in
R[M]_{[\alpha_n,+\infty)},\\
&e_{ij}(-\gamma_1)a_{ji}(\beta)=a_{ji}(\beta)+D_2,\quad
D_2\in\D_{>0}.
\end{align*}
We can write $\delta_1=\gamma_2+a_2+b_2$ for some
\begin{align*}
\gamma_2\in R[M]_{[\alpha_n,\epsilon_1)}\setminus
\B'(-\epsilon,l),\quad a_2\in R[M]_{[\epsilon_1,+\infty)},\quad
b_2\in\B'(-\epsilon,l).
\end{align*}

\smallskip If $\gamma_2=0$ then we are done because
\begin{align*}
&e_{ij}(-\gamma_1)e_{ji}(-\alpha)(e_{ji}(\beta)+D)e_{ij}(\alpha)=\\
&\1+\left[(a_{ij}(a_1+a_2)+a_{ji}(\beta)\right]+a_{ij}(b_1+b_2)+\left[D+D_1+D_2\right].
\end{align*}

\smallskip Therefore, there is no loss of generality in assuming
that $\gamma_2\not=0$. Then we derive elements $\gamma_3$, $a_3$,
$b_3$, $\delta_2$ and a matrix $D_3$ out from $\gamma_2$, $a_2$,
$b_2$, $a_1$, $b_1$ and $D+D_1+D_2$ in the same way $a_2$, $b_2$,
$\gamma_2$, $\delta_1$ and $D_2$ were derived out from $a_1$,
$b_1$, $\gamma_1$ and $D+D_1$, etc.

\smallskip If we show that $\gamma_p=0$ for some $p\in\NN$ then
$$
e_{ij}\left(-\gamma\right)e_{ij}(-\alpha)\left(e_{ji}(\beta)+D\right)=
\1+\left[a_{ij}(\alpha')+a_{ji}(\beta)\right]+a_{ij}(\beta')+D'
$$
where
\begin{align*}
&\gamma=\sum_{k=1}^{p-1}\gamma_k\in
R[M|D,\alpha,\beta]_{[\alpha_n,\epsilon_1)}\setminus
B'(-\epsilon,l),\\
&\alpha'=\sum_{k=1}^pa_k\in R[M]_{[\epsilon_1,+\infty)},\quad
\beta'=\sum_{k=1}^pb_k\in\B'(-\epsilon,l),\quad
D'=\sum_{k=1}^pD_k\in \D,
\end{align*}
and the lemma is proved.

\

Assume to the contrary that $\gamma_p\not=0$ for all $p\in\NN$. On
the other hand it follows from the definition of the elements
$\gamma_p$ that every element of $\supp(\gamma_{p+1})$ is
\emph{strictly} divisible in $M$ by a some element of
$\supp(\gamma_p)$. (In fact, we have
$\supp(D),\supp(D_1),\ldots\subset RM^+$ for all $p\ge1$.) Since
$M$ is an affine positive monoid, $\|\gamma_p\|\to\infty$ as
$p\to\infty$. But we also have $\gamma_k\in
R[M]_{[\alpha_n,+\infty)}$. Therefore, if $p$ is big enough, then
the radial direction of the support terms of $\gamma_p$ are almost
parallel to $\RR^{n-1}\oplus0\subset\RR^n$ and, in particular,
belong to $\B'(-\epsilon,l)$.

\smallskip The claim that the support monomials of $\gamma$, $A$,
$B$ and $D'$ are products of those of $\alpha$, $\beta$ and $D$ is
a consequence of the process of constructing these objects.
\end{proof}

\begin{lemma}\label{Com2}
Let $\epsilon_1,\epsilon,l$ be positive real numbers, $A\in\A(0)$
and $B\in\B(-\epsilon,l)$. Then
$$
E(\1+A+B)=\1+A'+B'+D'
$$
for some $A'\in\A(\epsilon_1)$, $B'\in\B(-\epsilon,l)$, $D'\in\D$
and $E\in\E_r(R[M])$ with  a representation $\bar E$ such that
$0\le\bar E_n<\epsilon_1$. Moreover, the support monomials of
$A'$, $B'$, $D'$ and of the factors in $\bar E$ are products of
the support monomials of $A$.
\end{lemma}

\begin{proof}
Let $A=(\alpha_{ij})$. For every pair of indices $i\not=j$ we let
$\bar\alpha_{ij}$ be the sum of those terms in the canonical
expansion of $\alpha_{ij}$ that have the $n$th coordinate
$<\epsilon_1$ and whose length is $<l$. We have a representation
of the form
$$
\left(\prod_{i\not=j}e_{ij}(-\bar\alpha_{ij})\right)(\1+A+B)=\1+A_1+B_1
$$
where:

\smallskip$\centerdot$ the order of factors is chosen arbitrarily,

\smallskip$\centerdot$ $A_1\in\A(0)$,

\smallskip$\centerdot$ $B_1\in\B(-\epsilon,l)$,

\

The inequality (\ref{Acute}) in Section \ref{Notation} implies
that $\B(-\epsilon,l)$ is stable under the multiplication by
elementary matrices of the form $e_{ij}(\lambda)$ with
$0\le\lambda_n$. Therefore, we can repeat the process with respect
to the matrix $\1+A_1+B_1$ etc. The standard elementary matrices
that are produced in this process are of the from
$e_{ij}(\lambda)$ with $0\le\lambda_n<\epsilon_1$. After $p$ steps
we will have a representation of the form
$$
E_p(\1+A+B)=\1+A_p+B_p
$$
where:

\smallskip$\centerdot$ $A_p\in\A(0)$ and $B_p\in\B(-\epsilon,l)$,

\smallskip$\centerdot$ $E_p\in\E_r(R[M])$, having a representation
$\bar E_p$ with $0\le(\bar E_p)_n<\epsilon_1$.

\smallskip$\centerdot$
if a support monomial of some non-diagonal entry of $A_p$ has the
$n$-th coordinate $<\epsilon_1$ and the length $<l$ then it is a
product of $p$ elements (maybe with repetitions) of $M^+$.

\

Because $M$ is affine positive, the lengths of the products
mentioned in the last condition above go to $\infty$ as
$p\to\infty$. In other words, if $p$ is big enough then the
mentioned support terms simply do not exist. That is, for $p$
large enough $\1+A_p+B_p=\1+A'+B'+D'$ for some
$A'\in\A(\epsilon)$, $B'\in\B(-\epsilon,l)$ and $D'\in\D$.

\smallskip As in the previous lemma, the claim that the support
monomials of $A'$, $B'$, $D'$ and of the factors in $\bar E$ are
products of the support monomials of $A$ and $B$ is a consequence
of the process by which these matrices have been constructed.
\end{proof}

To formulate the next result we introduce certain function
$\mathfrak l:\RR^3_{>0}\to\RR_{>0}$, where $\RR_{>0}$ is the set
of positive reals. For a triple
$(\epsilon_1,\epsilon_2,\epsilon)\in\RR_{>0}$ there exists a real
number $l(\epsilon_1,\epsilon_2,\epsilon)>0$ such that the
following implication holds:
\begin{equation}\label{large}
\begin{aligned} l\ge l(\epsilon_1,\epsilon_2,\epsilon),\
&A_1,A_2\in\A(-\epsilon_1),\
B\in\B(-\epsilon_2,l)\ \Longrightarrow\\
&A_1B,\ BA_2,\ A_1BA_2\in\B(-\epsilon_2-\epsilon,l).
\end{aligned}
\end{equation}
In fact, if $m_1\in\supp(A_1)$, $m_2\in\supp(A_2)$ and
$x\in\supp(B)$ then the inequality (\ref{Acute}) in Section
\ref{Notation} implies $|m_1x|,|m_2x|,|m_1m_2x|\ge l$. On the
other hand, none of the numbers $\Phi(m_1x)_n$, $\Phi(m_2x)_n$ and
$\Phi(m_1m_2x)_n$ can be less than $\Phi(-2\epsilon_1e_n+x)_n$
(switching do additive notation). Now if $l\gg0$, depending on
$\epsilon_1$, $\epsilon_2$ on $\epsilon$, then
$\Phi(-2\epsilon_1e_n+x)_n$ cannot be less than
$\Phi(x)_n-\epsilon$.

\smallskip The function $\mathfrak l$ is defined by
$(\epsilon_1,\epsilon_2,\epsilon)\mapsto
l(\epsilon_1,\epsilon_2,\epsilon)$.

\begin{proposition}\label{Com3}
Let:

\smallskip$\centerdot$ $\epsilon_1,\epsilon_2,\epsilon,l$ be positive real numbers with $l\ge\mathfrak
l(\epsilon_1,\epsilon_2,\epsilon)$,

\smallskip$\centerdot$
$i\not=j$ be natural numbers,

\smallskip$\centerdot$ $\alpha\in R[M]$ be a nonzero term with $|\alpha_n|<\epsilon_1$,

\smallskip$\centerdot$ $A\in\A(\epsilon_1)$, $B\in\B(-\epsilon_2,l)$ and
$D\in\D$.

\

\noindent Then:
$$
(\1+A+B+D) e_{ij}(\alpha)=e_{ij}(\alpha) E (\1+A_1+B_1+D_1)
$$
for some $A_1\in\A(\epsilon_1)$,
$B_1\in\B(-\epsilon_2-\epsilon,l)$, $D_1\in\D$ and
$E\in\E_r(R[M])$, having a representation $\bar E$ such that
$\min(\alpha_n,0)\le\bar E_n<\epsilon_1$. Moreover, the support
monomials of $A_1$, $B_1$, $D_1$ and of the factors in $\bar E$
are products of the support monomials of $\alpha$, $A$, $B$ and
$D$.
\end{proposition}

\smallskip\noindent(Observe, we do not exclude the case $\alpha\in R$.)

\begin{proof} Let $\beta$ be the $ji$-entry of $A$. Then $|\alpha_n|<\epsilon_1\le\beta_n$
and by Lemma \ref{Com1} we have a representation of the form
\begin{equation}\label{Com3a1}
\left(e_{ji}\left(\beta\right)+D\right)e_{ij}(\alpha)=e_{ij}(\alpha+\gamma)
\left(\1+A'+B'+D'\right)
\end{equation}
where $\gamma\in R[M]_{[\alpha_n,\epsilon_1)}$,
$A'\in\A(\epsilon_1)$, $B'\in\B(-\epsilon_2,l)$ and $D'\in\D$.

\smallskip We have
\begin{equation}\label{Com3a3}
A''=e_{ij}(\alpha-\gamma)\left(A-a_{ji}(\beta)\right)e_{ij}(\alpha)\in\A(0)
\end{equation}
because
$$
\supp(A'')\subset\supp(A)\cup\{\alpha x\ |\
x\in\supp(A)\}\cup\{\gamma x |\ x\in\supp(A)\}.
$$

In view of the implication (\ref{large}), we also have
\begin{equation}\label{Com3a2}
B''=e_{ij}(-\alpha-\gamma)B
e_{ij}(\alpha)\in\B(-\epsilon_2-\epsilon,l).
\end{equation}

Using (\ref{Com3a1}) and the definition of the matrices $A''$ and
$B''$, we can write:
\begin{align*}
&e_{ij}(-\alpha-\gamma)(\1+A+B+D)e_{ij}(\alpha)=\\
&e_{ij}(-\alpha-\gamma)\left(e_{ji}(\beta)+D\right)e_{ij}(\alpha)+A''+B''=\\
&\1+\left(A'+A''\right)+\left(B'+B''\right)+D'.
\end{align*}
We have $A'+A''+D'\in\A(0)$ by (\ref{Com3a3}) and
$B'+B''\in\B(-\epsilon_2-\epsilon,l)$ by (\ref{Com3a2}). By Lemma
\ref{Com2} we get a representation of the form:
$$
E\big(\1+(A'+A''+D')+(B'+B'')\big)=\1+A_1+B_1+D_1
$$
where: $A_1\in\A(\epsilon_1),\quad
B_1\in\B(-\epsilon_2-\epsilon,l)$, $D_1\in\D$, and
$E\in\E_r(R[M])$, having a representation $\bar E$ such that
$0\le\bar E_n<\epsilon_1$.

We finally get the desired representation:
$$
E e_{ij}(-\alpha-\gamma)(\1+A+B+D)e_{ij}(\alpha)=\1+A_1+B_1+D_1,
$$
that is
$$
(\1+A+B+D)e_{ij}(\alpha)=e_{ij}(\alpha)\left(e_{ij}(\gamma)\cdot
E^{-1}\right) \left(\1+A_1+B_1+D_1\right).
$$

\smallskip That the support monomials of $A_1$, $B_1$,
$D_1$ and of the factors in $e_{ij}(\gamma)\cdot E^{-1}$ are
products of the support monomials of $\alpha$, $A$, $B$ and $D$
follows from the corresponding claims in Lemmas \ref{Com1} and
\ref{Com2} and the way these lemmas are used in the argument
above.
\end{proof}

\subsection{Almost separation}\label{lastsection} Finally, here we prove
Theorem \ref{elparation}.

In addition to the objects and the conditions on them, listed in
Section \ref{Notation}, we now require that $M$ is normal and
$\gp(M)=\ZZ^n$.

Also, we extend in the obvious way to the monoid ring
$R[(M_*)^{c^{-\infty}}]$ the terminology and notation that was
introduced in Section \ref{Notation} for $R[M]$.

Assume $\RR_+M=C_1\cup C_2$ where $C_1=\{z\in \RR_+M\ |\
z_n\le0\}$ and $C_2=\{z\in \RR_+M\ |\ z_n\ge0\}$.

\smallskip Fix a real number $\epsilon>0$. As in Theorem
\ref{elparation}, we let $M_1(\epsilon)=\RR_+M\cap
C_1(\epsilon)\cap M$ and $M_2(\epsilon)=\RR_+M\cap
C_2(\epsilon)\cap M$.

\smallskip Let $c$ be a natural number $\ge2$.

We want to prove the inclusion:
\begin{equation}\label{finally!}
\E_r(R[(M_*)^{c^{-\infty}}])\subset\E_r(R[(M_1(\epsilon)_*)^{c^{-\infty}}])
\SL_r(R[(M_2(\epsilon)_*)^{c^{-\infty}}]),
\end{equation}
the left hand side being considered in
$\SL_r(R[(M_*)^{c^{-\infty}}])$.

\begin{lemma}\label{goodrep}
For (\ref{finally!}) it is enough to consider the matrices
$E=\Pi_{k=1}^{s}e_{i_kj_k}(\alpha_k)$ where:
\begin{itemize}
\item[(a)] $\alpha_k$ are terms in $R[(M_*)^{c^{-\infty}}]$,
\item[(b)] $(\alpha_k)_n\in\ZZ$, \item[(c)] $(\alpha_k)_n<0\
\Longrightarrow\ (\alpha_k)_n=-1$, \item[(d)] $(\alpha_k)_n>0,\
\beta\in R[(M_*)^{c^{-\infty}}],\ (\alpha_k\beta)_n=1
\Longrightarrow\ \alpha_k\beta\in(M_1(\epsilon)_*)^{c^{-\infty}}$.
\end{itemize}
\end{lemma}

\begin{proof}
Consider any matrix
$E'=\prod_ke_{i_kj_k}(\alpha'_k)\in\E_r(R[(M_*)^{c^{-\infty}}])$.
In view of the 1st Steinberg relation (Section \ref{Kbaground}) we
can assume that $\alpha'_k\in R[(M_*)^{c^{-\infty}}]$ are terms.
Assume $\alpha_k'=a_k\mu_k$ for some $a_k\in R$ and $\mu_k\in
(M_*)^{c^{-\infty}}$. It is enough to consider the matrix
$(c^j)_*(E')$ for some $j\gg0$. Therefore, there is no loss of
generality also in assuming that $\mu_k\in M_*$ for all $k$.
Moreover, by taking $j$ sufficiently large we can make the lengths
$\|\mu_k\|$ large enough so that the condition (d) is satisfied.
In more detail, we have $0\ll\|\alpha'_k\|\le\|\alpha'_k\beta\|$
for any monomial $\beta\in R[(M_*)^{c^{-\infty}}]$, the second
inequality being implied by (\ref{Acute}) in Section
\ref{Notation}. But a long monomial with the $n$th coordinate $=1$
must be almost parallel to the hyperplane $\mathcal
H=\RR^{n-1}\oplus0$, or equivalently, must belong to the submonoid
$(M_1(\epsilon)_*)^{c^{-\infty}}\subset(M_*)^{c^{-\infty}}$.

At this point we have reached the situation when all but the
condition (c) are satisfied. Now the mentioned condition is taken
care of as follows.

The normality of $M$ and the equality $\gp(M)=\ZZ^n$
(equivalently, the condition $M=\RR_+M\cap\ZZ^n$) imply the
surjectivity of the monoid homomorphism $M_*\to\ZZ$,
$\mu\mapsto\mu_n$. Therefore, by Lemma \ref{decomposition} for
every $\mu_k$ with $(\mu_k)_n<0$ there exists a decomposition of
the form (in additive notation):
$$
\mu_k=\sum_i\mu_{ki},\quad\mu_{ki}\in (M_*)^{c^{-\infty}}\cap
h^{-1}(-1).
$$
Using the 3rd Steinberg relation (Section \ref{Kbaground}) the
matrices $e_{i_kj_k}(\alpha'_k)$ with $(\alpha_k')_n<0$ can
correspondingly be represented as products of matrices of the form
$$
e_{pq}(a_k),\ e_{pq}(\mu_{k1}),\ e_{pq}(\mu_{k2}),\ldots
$$
Substituting in the product $\prod_ke_{i_kj_k}(\alpha'_k)$ these
representations correspondingly for the factors
$e_{i_kj_k}(\alpha'_k)$, $(\alpha_k')_n<0$, we arrive at the
desired representation.
\end{proof}

\begin{proof}[Proof of the equality (\ref{finally!})]
Products of elementary matrices of the form mentioned in Lemma
\ref{goodrep} will be called \emph{admissible representations}.

Let $E\in\E_r(R[(M_*)^{c^{-\infty}}]$, having an admissible
representation $\bar E=\prod_{k=1}^{s}e_{i_kj_k}(\alpha_k)$. We
want to show
\begin{equation}\label{GOOD}
E\in\E_r(R[(M_1(\epsilon)_*)^{c^{-\infty}}])
\SL_r(R[(M_2(\epsilon)_*)^{c^{-\infty}}]).
\end{equation}

Let $M'\subset(M_*)^{c^{-\infty}}$ be the submonoid generated by
$\cup_k\supp(\alpha_k)$ and $\tilde M\subset(M_*)^{c^{-\infty}}$
be the submonoid generated by $M\cup M'$.

It is important that the elements of $\tilde M$ have
\emph{integral} $n$th coordinate.

An admissible representation of $E$ whose factors have support
monomials in $M'$ will be called \emph{good}.

Assume $(\alpha_k)_n\le a$ for some $a\ge0$. Let
$$
\alpha_{k_1},\ldots,\alpha_{k_p},\qquad1\le k_1<k_2<\cdots<k_p\le
s,
$$
be determined by the condition:
$$
\left(\alpha_{k_1}\right)_n,\ldots,\left(\alpha_{k_p}\right)_n=a.
$$
In this situation we say that the representation $\bar E$ is
\emph{$(a,p)$-bounded}.

\smallskip Consider the lexicographic order on
$\ZZ_+\times\ZZ_+$. For any pair $(a',p')$ with $(a,p)\le(a',p')$
we also say that $\bar E$ is $(a',p')$-bounded.

\smallskip The proof is by induction on the bounding pairs.

\smallskip If $a=1$ then (\ref{GOOD}) follows from the condition
(d) in Lemma \ref{goodrep}: in this situation
$E\in\E_r(R[(M_1(\epsilon)_*)^{c^{-\infty}}])$.

\smallskip So we can assume $a\ge2$ and that
$$
E\in\E_r(R[(M_1(\epsilon)_*)^{c^{-\infty}}])
\SL_r(R[(M_2(\epsilon)_*)^{c^{-\infty}}])
$$
whenever $E$ has an $(a',p')$-bounded good representation for some
$(a',p')<(a,p)$.

\smallskip It is enough to prove the existence of a representation
of the form:
\begin{equation}\label{E=YZ}
\begin{aligned}
E=YZ,\quad &Y\in\E_r(R[(M_1(\epsilon)_*)^{c^{-\infty}}],\
\text{having an}\ (a',p')\text{-bounded}\\
&\text{good representation}\ \bar Y\ \text{for some}\ (a',p')<(a,p),\ \text{and}\\
&Z\in\SL_r(R[(M_2(\epsilon)_*)^{c^{-\infty}}]).
\end{aligned}
\end{equation}

\

There is no loss of generality in assuming that $k_p<s$ for
otherwise
$$
E=\left(Ee_{i_sj_s}(-\alpha_s)\right)e_{i_sj_s}(\alpha_s)
$$
and $Ee_{i_sj_s}(-\alpha_s)$ obviously has an $(a',p')$-bounded
good representation for some $(a',p')<(a,p)$.

\smallskip Fix positive real numbers $\epsilon_2$ and $\epsilon'$
so that $\epsilon_2+(s-k_p)\epsilon'=\epsilon$. Also, fix a real
number $l>0$, sufficiently large with respect to the numbers
$$
a,\ \epsilon_2,\ \epsilon',\ \epsilon_2+\epsilon',\
\epsilon_2+2\epsilon',\ \ldots,\ \epsilon_2+(s-k_p-1)\epsilon'.
$$
We apply Proposition \ref{Com3} to the product
$$
e_{i_{k_p}j_{k_p}}(\alpha_{k_p})e_{i_{k_p}+1j_{k_p}+1}(\alpha_{k_p+1})
$$
where in the notation of Proposition \ref{Com3}:

\smallskip$\centerdot$ the r\^ole of $M$ is played by $\tilde M$,

\smallskip$\centerdot$ $\epsilon_1=a$, $\epsilon_2=\epsilon_2$ and
$\epsilon=\epsilon'$,

\smallskip$\centerdot$ $\1+A+B+D=\1+A+0+0=e_{i_{k_p}j_{k_p}}(\alpha_{k_p})$,

\smallskip$\centerdot$ $e_{ij}(\alpha)=e_{i_{k_p}+1j_{k_p}+1}(\alpha_{k_p+1})$.

\

We get
\begin{align*}
e_{i_{k_p}j_{k_p}}(\alpha_{k_p})e_{i_{k_p}+1j_{k_p}+1}(\alpha_{k_p+1})=
e_{i_{k_p}+1j_{k_p}+1}(\alpha_{k_p+1})E_1(\1+A_1+B_1+D_1)
\end{align*}
for some $A_1\in\A(a)$, $B_1\in\B(-\epsilon_2-\epsilon',l)$,
$D_1\in\D$, and $E_1\in\E_r(R[\tilde M])$, having a good
representation $\bar E_1$ with $(\bar E_1)_n<a$.

\smallskip Using Proposition \ref{Com3}, we can find inductively matrices
\begin{align*}
A_t\in\A(a),\quad B_t\in\B(-\epsilon_2-t\epsilon'&,l),\quad
D_t\in\D,\\
&t\in\{1,\ldots,s-k_p-1\},
\end{align*}
starting with the triple $A_1,B_1,D_1$ above, so that the
following holds for each $t$:
\begin{align*}
(\1+A_t+B_t+&D_t)e_{i_{k_p}+tj_{k_p}+t}(\alpha_{k_p+t})=\\
&e_{i_{k_p}+tj_{k_p}+t}(\alpha_{k_p+t})E_{t+1}(\1+A_{t+1}+B_{t+1}+D_{t+1}),
\end{align*}
where $A_{t+1}\in\A(a)$,
$B_{t+1}\in\B(-\epsilon_2-(t+1)\epsilon',l)$, $D_{t+1}\in\D$, and
$E_{t+1}\in\E_r(R[\tilde M])$, having a good representation $\bar
E_{t+1}$ with $(\bar E_{t+1})_n<a$.

\smallskip We have
\begin{align*}
e_{i_{k_p}j_{k_p}}(\alpha_{k_p})\prod_{t=i_{k_p}+t}^se_{i_tj_t}(\alpha_t)=
\mathcal E(\1+A_s+B_s+D_s)
\end{align*}
for some $\mathcal E\in\E_r(R[\tilde M])$ having a good
representation $\bar{\mathcal E}$ with $\bar{\mathcal E}_n<a$.
Hence a representation $E=YZ$ of the form (\ref{E=YZ}) where:

\smallskip$\centerdot$ $Y=\left(\prod_{t=1}^{k_p-1} e_t(\alpha_t)\right)\mathcal
E$,

\smallskip$\centerdot$ $Z=\1+A_s+B_s+D_s$.
\end{proof}

\end{document}